\documentclass[leqno,12pt]{article}
\usepackage{latexsym,amsmath,amsfonts,amssymb}
\usepackage[dvips]{color}

\scrollmode

\newtheorem{thm}{Theorem}[section]
\newtheorem{lem}[thm]{Lemma}
\newtheorem{cor}[thm]{Corollary}
\newtheorem{rem}[thm]{Remark}
\newtheorem{rems}[thm]{Remarks}
\newtheorem{prop}[thm]{Proposition}
\newtheorem{df}[thm]{Definition}
\newtheorem{dfs}[thm]{Definitions}
\newtheorem{ex}[thm]{Example}

\newcommand{\wh}{\widehat}
\newcommand{\wt}{\widetilde}
\newcommand{\vS}{\varSigma}
\newcommand{\vO}{\varOmega}
\newcommand{\vT}{\varTheta}
\newcommand{\ov}{\overline}
\newcommand{\mf}{\mathfrak}
\newcommand{\vY}{\varUpsilon}
\newcommand{\vX}{{\varPsi}}
\def\N{{\mathbb N}}
\def\R{{\mathbb R}}
\def\Q{{\mathbb Q}}

\newcommand{\cnp}{\{N_t\}_{t\in\R_+}}
\newcommand{\cnpw}{counting process }

\newcommand{\clap}{\{T_n\}_{n\in\N_0}}

\newcommand{\cip}{\{W_n\}_{n\in\N}}

\newcommand{\E}{\mathbb{E}}

\newcommand{\bdot}{\bullet}
\newcommand{\leb}{\lambda}

\definecolor{greych}{rgb}{0.9,0.9,0.9}

\newcommand{\auth}[1]{\textcolor{black}{\texttt{#1}}}
\newcommand{\artp}[1]{\textrm{#1}}
\newcommand{\artn}[1]{\textbf{#1}}
\newcommand{\artj}[1]{\textsl{#1}}
\newcommand{\artb}[1]{\emph{#1}}

\usepackage[pdftex,colorlinks={false}]{hyperref}

\title{SOME CHARACTERIZATIONS OF\\MIXED RENEWAL PROCESSES}

\author{\auth{D.P. Lyberopoulos}\thanks{The author is indebted to the Public Benefit Foundation \textsc{Alexander S. Onassis}, which supported this research, under the Programme 
of Scholarships for Hellenes.}\;\;{\normalsize and}\;\;\auth{N.D. Macheras}}

\date{{\small\textsf{\today}}}
\begin{document}
\maketitle

\marginpar{lm6z15}

\begin{abstract}
\noindent Some characterizations of mixed renewal processes in terms of exchangeability and of different types of disintegrations are given. As a consequence, an existence result for mixed renewal processes, providing also a new construction for them, is obtained. As an application, some concrete examples of constructing such processes are presented and the corresponding disintegrating measures are explicitly computed.\smallskip

\smallskip

\end{abstract}

\section{Introduction}\label{intro}

Mixed renewal processes (MRPs for short) may serve as a source of challenging theoretical problems, since they are generalizations of mixed Poisson processes (MPPs for short) and closely connected with the concept of exchangeable stochastic processes (cf. e.g. \cite{hu}), as well as a useful tool for modelling real life situations, such as those emerging in actuarial practice (cf. e.g. \cite{se}, pages 164-165).

In Section \ref{char} we introduce a new definition of MRPs (see Definition \ref{rd}) being in line with that of MPPs with parameter $\vT$. Such a definition seems to be a proper one as it involves explicitly the structural parameter $\vT$, which is usually essential in the study of risk-theoretical problems. Since conditioning is involved in this definition of MRPs, it seems to be natural to ask about the structural role of disintegrations in this field. For this reason, we recall the definitions of different types of disintegrations (see Definitions \ref{rcp}) and provide some characterizations of MRPs via disintegrations (see Proposition \ref{31} and Corollary \ref{25p}). By means of these disintegration results we can reduce MRPs to ordinary renewal processes for the disintegrating probability measures, showing in this way that Definition \ref{rd} is the natural one for MRPs. 
Furthermore, we prove a characterization of conditionally independent and conditionally independent and identically
distributed (conditionally i.i.d. for short) families of measurable maps (see Proposition \ref{nid0}), implying an extension (see Theorem \ref{19}) of a basic result, i.e. of Proposition 4.4 of our previous work \cite{lm1v}.

The second definition of MRPs investigated in this paper is due to Huang \cite{hu}, see Definition \ref{hua}. 
In Section \ref{exch}, we first obtain, by an application of Proposition \ref{nid0}, some characterizations of exchangeability in terms of different types of disintegrations, see Theorem \ref{21}. As a consequence, some further characterizations of MRPs in terms of exchangeability and of disintegrations are deduced 
(see Theorem \ref{9}). Theorem \ref{9} provides amongst others a detailed discussion of the relation between the two definitions of MRPs and shows that in most cases appearing in applications both definitions coincide.

In Section \ref{con} we first give some examples to show that some of the assumptions of Theorems \ref{21} and \ref{9} are essential for the validity of all the equivalences obtained therein (see Examples \ref{exa} to \ref{exd}). Next, we provide a construction of non-trivial probability spaces admitting MRPs (see Example \ref{12}), extending a similar construction for MPPs, see \cite{lm5jmaslo}, Theorem 3.1. As an application, we give concrete examples of MRPs satisfying all assumptions of Theorem \ref{9}, and we compute the corresponding disintegrating measures explicitly.

\section{Preliminaries}\label{prel}

By $\N$ is denoted the set of all natural numbers and $\N_0:=\N\cup\{0\}$.  The symbol $\R$ stands for the set of all real numbers, while $\R_+:=\{x\in\mathbb{R}:x\geq0\}$. If $d\in\N$, then $\R^d$ denotes the Euclidean space of dimension $d$.

Given a probability space $(\vO,\vS,P)$, a set $N\in\vS$ with $P(N)=0$ is called a $P$-{\bf null set} (or a null set for simplicity). The family of all $P$-null sets is denoted by $\vS_0$. 
For any two $\vS$-$T$-measurable maps $X,Y:\vO\longrightarrow\vY$ we write $X=Y\;P-\mbox{a.s.}$, if $\{X\not=Y\}\in\vS_0$.
 
If $A\subseteq\vO$, then $A^c:=\vO\setminus A$, while $\chi_A$ denotes the indicator (or characteristic) function of the set $A$. The identity map from $\vO$ onto itself is denoted by ${\mathrm{id}}_{\vO}$. The $\sigma$-algebra generated by a family $\mathcal{G}$ of subsets of $\vO$ is denoted by $\sigma(\mathcal{G})$. A $\sigma$-algebra $\mathcal{A}$ is {\bf countably generated} if there exists a countable family $\mathcal{G}$ of subsets of $\vO$ such that $\mathcal{A}=\sigma(\mathcal{G})$.

For any Hausdorff topology $\mf{T}$ on $\vO$ by ${\mf B}(\vO)$ is denoted the 
{\bf Borel $\sigma$-algebra} on $\vO$, i.e. the $\sigma$-algebra generated by $\mf{T}$. 
By ${\mf B}:={\mf B}(\R)$, ${\mf B}_d:={\mf B}(\R^d)$ and $\mf{B}_{\N}:={\mf B}(\R^{\N})$ is denoted the Borel $\sigma$-algebra of subsets of $\R$, $\R^d$ and $\R^{\N}$, respectively, while ${\mathcal{L}}^{1}(P)$ stands for the family of all real-valued $P$-integrable functions on $\vO$. Functions that are $P$-a.s. equal are not identified.

A $\vS$-$\mf{B}$-measurable function $X$ from $\vO$ into $\R$ is called a {\bf random variable}, while if $X$ is a $\vS$-$\mf{B}(\ov\R)$-measurable map from $\vO$ into $\ov\R$, then it is said to be an {\bf extended random variable}. Also recall that any $\vS$-$\mf{B}_d$-measurable function $\vT$ from $\vO$ into $\R^d$ is called a 
{\bf ($d$-dimensional) random vector}.

The probability measure $P$ is said to be {\bf perfect} if for any random variable $X$ on $\vO$ there exists a Borel set $B\subseteq X(\vO):=\{X(\omega):\;\omega\in\vO\}$ such that $P(X^{-1}(B))=1$. 

Given two probability spaces $(\vO,\vS,P)$ and $(\vY,T,Q)$ as well as a $\vS$-$T$-measurable map $X:\vO\longrightarrow\vY$ we denote by $\sigma(X):=\{X^{-1}(B): B\in T\}$ the $\sigma$-algebra generated by $X$, while 
$\sigma(\{X_i\}_{i\in I}):=\sigma\bigl(\bigcup_{i\in{I}}\sigma(X_i)\bigr)$ stands for the $\sigma$-algebra generated by a family $\{X_i\}_{i\in I}$ of $\vS$-$T$-measurable maps from $\vO$ into $\vY$.

Setting $T_X:=\{B\subseteq\vY:\;X^{-1}(B)\in\vS\}$ for any given $\vS$-$T$-measurable map $X$ from $\vO$ into $\vY$, we clearly get that $T\subseteq T_X$. Denote by 
$P_X:T_X\longrightarrow\R$ the image measure of $P$ under $X$. The restriction of $P_X$ to $T$ is denoted again by $P_X$. By $\mathbf{K}(\theta)$ is denoted an arbitrary probability distribution on $\mf{B}$ with parameter $\theta\in\vX$. In particular, $\mathbf{P}(\theta)$, $\mathbf{Exp}(\theta)$ and $\mathbf{Ga}(\gamma,\alpha)$, where $\theta,\gamma,\alpha$ are positive parameters, stand for the law of Poisson, exponential and gamma distribution, respectively (cf. e.g. \cite{sch}).

If $X\in{\mathcal{L}}^{1}(P)$ and $\mathcal{F}$ is a $\sigma$-subalgebra of $\vS$, then each function $Y\in{\mathcal L}^1(P\mid\mathcal{F})$ satisfying for each $A\in{\mathcal{F}}$ the equality $\int_AX\,dP=\int_AY\,dP$ is said to be a {\bf version of the conditional expectation} of $X$ with respect to (or given) $\mathcal{F}$ and it will be denoted by $\E_{P}[X\mid\mathcal{F}]$. For $X:=\chi_B\in\mathcal{L}^1(P)$ with $B\in\vS$ we set $P(B\mid\mathcal{F}):=\E_P[\chi_B\mid\mathcal{F}]$.

By $(\vO\times\vY,\vS\otimes{T},P\otimes{Q})$ is denoted the product probability space of $(\vO,\vS,P)$ and $(\vY,T,Q)$, and by $\pi_{\vO}$ and $\pi_{\vY}$ the canonical projections from $\vO\times\vY$ onto $\vO$ and $\vY$, respectively.

Given two measurable spaces $(\vO,\vS)$ and $(\vY,T)$, a {\bf $T$-$\vS$-Markov kernel} is a function $k$ from $T\times\vO$ into $\R$ satisfying the following conditions:
\begin{description}
\item[(k1)] The set-function $k(\cdot,\omega)$ is a probability measure on $T$ for any fixed $\omega\in\vO$.
\item[(k2)] The function $\omega\longmapsto k(B,\omega)$ is $\vS$-measurable for any fixed $B\in{T}$.
\end{description}
Let be given a $\vS$-$T$-measurable map $X$ from $\vO$ into $\vY$ and a $\sigma$-subalgebra 
$\mathcal{F}$ of $\vS$. A {\bf conditional distribution of $X$ over $\mathcal{F}$} is a $T$-$\mathcal{F}$-Markov kernel $k$ satisfying for each $B\in T$ condition
$$
k(B,\cdot)=P(X^{-1}(B)\mid\mathcal{F})(\cdot)\quad{P}\mid\mathcal{F}-\mbox{a.s.}.
$$
Such a Markov kernel $k$ will be denoted by $P_{X\mid\mathcal{F}}$. 
In particular, if $(\vX,Z)$ is a measurable space, $\vT$ is a $\vS$-$Z$-measurable map from $\vO$ into $\vX$ and $\mathcal{F}:=\sigma(\vT)$, then the function $P_{X\mid\vT}:=P_{X\mid\sigma(\vT)}$ is called a {\bf conditional distribution of $X$ given $\vT$}. Note that if $\vY$ is a {\bf Polish space} (i.e. a topological space homeomorphic to a complete separable metric space) then a conditional distribution of $X$ over $\mathcal{F}$ always exists (cf. e.g. \cite{du}, Theorem 10.2.2).

Clearly, for every $T$-$Z$-Markov kernel $k$, the map $K(\vT)$ from $T\times\vO$ into $\R$ defined by means of
$$
K(\vT)(B,\omega):=(k(B,\cdot)\circ\vT)(\omega)
\quad\mbox{for any}\;\;B\in T\;\;\mbox{and}\;\;\omega\in\vO 
$$
is a $T$-$\sigma(\vT)$-Markov kernel. In particular, for $(\vY,T)=(\R,\mf{B})$ its associated probability measures $k(\cdot,\theta)$ for $\theta=\vT(\omega)$ with $\omega\in\vO$ are distributions on $\mf{B}$ and so we may write $\mathbf{K}(\theta)(\cdot)$ instead of $k(\cdot,\theta)$. 
Consequently, in this case $K(\vT)$ will be denoted by $\mathbf{K}(\vT)$.

For any $\sigma$-subalgebra $\mathcal{F}$ of $\vS$, we say that two $T$-${\mathcal{F}}$-Markov kernels $k_i$, for $i\in\{1,2\}$, are {\bf $P{\mid\mathcal{F}}$-equivalent} and we write $k_1=k_2$ $P{\mid\mathcal{F}}$-a.s., if there exists a $P$-null set $N\in{\mathcal{F}}$ such that for any $\omega\notin N$ and $B\in T$ the equality $k_1(B,\omega)=k_2(B,\omega)$ holds true. 

From now on $(\vO,\vS,P)$ is a probability space, while $(\vY,T)$ and $(\vX,Z)$ are measurable spaces, all of them arbitrary but fixed.

\section{Characterizations of mixed renewal processes via disintegrations}\label{char}

Let $\clap$ be a family of random variables, and let $\vO_T$ be a $P$-null set such that for each $\omega\in\vO\setminus\vO_T$ the sequence $\{T_n(\omega)\}_{n\in\N}$ is strictly increasing and $T_0(\omega)=0$. Then the family $\cnp$ defined by $N_t:=\sum_{n\in\N}\chi_{\{T_n\leq t\}}$ for each $t\in\R_+$ is the counting process {\em induced by $\clap$}. 

A family $\cnp$ is called a {\bf counting process} with exceptional $P$-null set $\vO_N$ if outside $\vO_N$ it takes values in $\N_0\cup\{\infty\}$, has right-continuous paths, presents jumps of size (at most) one, vanishes at $t=0$ and increases to infinity. Also recall that if $\cnp$ is a counting process with exceptional null set $\vO_N$, 
then the sequence $\clap$ defined by 
$$
T_n:=\inf\{t\in\R_+:N_t=n\}\quad\mbox{for each}\quad n\in\N_0
$$
is as above with $\vO_T=\vO_N$, while the sequence $\cip$, given by $W_n:=T_n-T_{n-1}$ for each $n\in\N$, is called the {\bf interarrival process} induced by $\clap$. Obviously, the exceptional null set of $\cip$ satisfies $\vO_{W}=\vO_{T}$. 
\smallskip

Recall now that a family $\{\vS_i\}_{i\in{I}}$ of $\sigma$-subalgebras of $\vS$ is {\bf $P$-conditionally (stochastically) independent} over a $\sigma$-subalgebra $\mathcal{F}$ of $\vS$, if for each $n\in\N$ with $n\geq 2$ we have
$$
P(E_{1}\cap\cdots\cap E_{n}\mid\mathcal{F})=\prod_{j=1}^nP(E_{{j}}\mid\mathcal{F})
\quad P\mid\mathcal{F}-\mbox{a.s.}
$$
whenever $i_1,\ldots,i_n$ are distinct members of $I$ and $E_j\in\vS_{i_j}$ for every 
$j\leq n$.
 
Let be given a family $\{X_i\}_{i\in I}$ of $\vS$-$T$-measurable maps from $\vO$ into $\vY$. We say that $\{X_i\}_{i\in I}$ is {\bf $P$-conditionally (stochastically) independent} over a $\sigma$-algebra $\mathcal{F}\subseteq\vS$, if the family $\{\sigma(X_i)\}_{i\in I}$ of $\sigma$-algebras is $P$-conditionally independent over $\mathcal{F}$. 

The family $\{X_i\}_{i\in I}$ is {\bf $P$-conditionally identically distributed} over $\mathcal{F}$, if 
$$
P\bigl(F\cap X_i^{-1}(B)\bigr)=P\bigl(F\cap X_j^{-1}(B)\bigr)
$$
whenever $i,j\in I$, $F\in\mathcal{F}$ and $B\in T$.

Furthermore, if $\vT$ is a $\vS$-$Z$-measurable map from $\vO$ into $\vX$, we say that 
$\{X_i\}_{i\in I}$ is {\bf $P$-conditionally (stochastically) independent or identically distributed given} $\vT$, if it is conditionally independent or identically distributed over the $\sigma$-algebra $\sigma(\vT)$.

{\em Throughout what follows, unless it is stated otherwise, $\vT$ is a $\vS$-$Z$-measurable map from $\vO$ into $\vX$, and we simply write ``conditionally'' in the place of ``conditionally given $\vT$'' whenever conditioning refers to $\vT$.}
{\em Moreover, $\cnp$ is a counting process and without loss of generality we may and do assume that $\vO_N=\emptyset$.} 
\smallskip

The counting process $\cnp$ is said to be a {\bf $P$-renewal process} with interarrival time distribution $\mathbf{K}(\theta_0)$, where $\theta_0\in{\R^d}$ 
is a parameter (or a $(P,\mathbf{K}(\theta_0))$-RP for short), if its associated interarrival times $W_n$, $n\in\N$, are independent and $\mathbf{K}(\theta_0)$-distributed under the probability measure $P$.

\begin{df}\label{rd}
\normalfont
The counting process $\cnp$ is said to be a {\bf mixed renewal process} on $(\vO,\vS,P)$ with parameter the $d$-dimensional ($d\in\N$) random vector $\vT$ and interarrival time conditional distribution $\mathbf{K}(\vT)$ (or a $(P,\mathbf{K}(\vT))$-MRP for short), if $\cip$ is $P$-conditionally independent and 
$$
P_{W_n\mid\vT}=\mathbf{K}(\vT)\quad P\mid\sigma(\vT)-\mbox{a.s.}
$$
for all $n\in\N$. 
\smallskip

In particular, for $(\vX,Z)=(\R,\mf{B})$ and $P_{\vT}((0,\infty))=1$ a counting process $\cnp$ is a {\bf $P$-mixed Poisson process} on $(\vO,\vS,P)$ (or a $P$-MPP for short) with parameter the random variable $\vT$, if it has $P$-conditionally stationary independent increments (cf. e.g. \cite{sch}, Section 4.1, page 86 for the definition), such that
$$
P_{N_t\mid\vT}=\mathbf{P}(t\vT)\quad P\mid\sigma(\vT)-\mbox{a.s.}
$$ 
holds true for each $t\in(0,\infty)$. 
\end{df}
Note that, for $d=1$, $P_{\vT}((0,\infty))=1$ and $\mathbf{K}(\vT)=\mathbf{Exp}(\vT)$ $P\mid\sigma(\vT)$-a.s. the $(P,\mathbf{K}(\vT))$-MRP $\cnp$ becomes a $P$-MPP with parameter $\vT$ (see \cite{lm1v}, Proposition 4.5).

\begin{dfs}\label{rcp}
\normalfont
{\bf (a)} 
Let $Q$ be a probability measure on $T$. A family $\{P_y\}_{y\in\vY}$ of probability measures on $\vS$ is called a {\bf disintegration} of $P$ over $Q$ if
\begin{description}
\item[(d1)] for each $D\in\vS$ the map $y\longmapsto P_y(D)$ is $T$-measurable;
\item[(d2)] $\int P_{y}(D)Q(dy)=P(D)$ for each $D\in\vS$.
\end{description}
If $f:\vO\longrightarrow\vY$ is an inverse-measure-preserving function (i.e. $P(f^{-1}(B))=Q(B)$ for each $B\in{T}$), a disintegration $\{P_{y}\}_{y\in\vY}$ of $P$ over $Q$ is called {\bf consistent} with $f$ if, for each $B\in{T}$, the equality $P_{y}(f^{-1}(B))=1$ holds for $Q$-almost every $y\in B$.

\noindent{\bf (b)}
Let $\mathcal{F}$ be a $\sigma$-subalgebra of $\vS$ and $R:=P\mid\mathcal{F}$. 
A {\bf subfield r.c.p.} for $P$ over $R$ (see \cite{fa}, Section 2) is a family $\{P_{\omega}\}_{\omega\in\vO}$ of probability measures on $\vS$ satisfying the following conditions: 
\begin{description}
\item[(sf1)] for each $E\in\vS$ the map $\omega\longmapsto P_{\omega}(E)$ is $\mathcal{F}$-measurable;
\item[(sf2)] $\int_FP_{\omega}(E)R(d\omega)=P(E\cap F)$ for all $F\in\mathcal{F}$ and $E\in\vS$.
\end{description}
\end{dfs}

\begin{rems}\label{mag}
\normalfont
{\bf (a)}
If $\vS$ is countably generated and $P$ is perfect, then there always exists a disintegration $\{P_{y}\}_{y\in\vY}$ of $P$ over $Q$ consistent with any inverse-measure-preserving map $f$ from $\vO$ into $\vY$ providing that $T$ is countably generated (see \cite{fa}, Theorems 6 and 3) and a subfield r.c.p. (see \cite{fa}, Theorems 6 and 2). So, in most cases appearing in applications (e.g. Polish spaces) these two types of disintegrations always exist.

\noindent{\bf (b)}
Let $\{P_y\}_{y\in\vY}$ be a disintegration of $P$ over $Q$, and let $f$ be an inverse-measure-preserving map from $\vO$ into $\vY$. Then the following are equivalent:
\begin{eqnarray}
\lefteqn{\quad\;\;\{P_y\}_{y\in\vY}\;\;\mbox{is consistent with}\;\;f}\label{mag1}\\
&&P(A\cap f^{-1}(B))=\int_B P_{y}(A)Q(dy)
\quad\mbox{for each}\;\;A\in\vS\;\mbox{and}\;\;B\in T\label{mag2}\\
&&\mbox{For each}\;\;A\in\vS\quad
\E_P[\chi_A\mid\sigma(f)]=P_{\bdot}(A)\circ f\quad P\mid\sigma(f)-\mbox{a.s.}.\label{mag3}
\end{eqnarray}
\noindent{\bf{(c)}}
Let $X$ be a $\vS$-$T$-measurable map from $\vO$ into $\vY$, let $\{P_{\theta}\}_{\theta\in\vX}$ be a disintegration of $P$ over $P_{\vT}$ consistent with $\vT$, and let $k$ be a $T$-$Z$-Markov kernel. If $k(\cdot,\theta)$ is the distribution of $X$ under $P_{\theta}$ for $\theta\in\vX$, then the map $K(\vT)$ is a conditional distribution of $X$ given $\vT$, since by condition (\ref{mag3}) of (b) we get for $A=X^{-1}(B)$ with $B\in T$ that $P_{X\mid\vT}(B,\cdot)=K(\vT)(B,\cdot)\;$ $P\mid\sigma(\vT)$-a.s..

\noindent{\bf{(d)}}
Conversely, in the special case where $\vS$ is countably generated and $(\vY,T)=(\R,\mf{B})$, given $\{P_{\theta}\}_{\theta\in\vX}$ as in {(c)}, we get that for each conditional distribution $\mathbf{K}(\vT)$ of $X$ given $\vT$, there exists an {\em essentially unique} probability distribution $(P_{\theta})_X$ of $X$, for $\theta\in\vX$, such that for each $B\in{\mf{B}}$ we have
$$
\mathbf{K}(\vT)(B,\cdot)=(P_{\bdot})_X(B)\circ\vT\quad P\mid\sigma(\vT)-\mbox{a.s.}.
$$
In fact, by applying a monotone class argument, it can be easily seen that the disintegration is essentially unique in the sense that if $\{P_{\theta}^{\prime}\}_{\theta\in\vX}$ is any other disintegration of $P$ over $P_{\vT}$ which is consistent with $\vT$, then $P_{\theta}=P_{\theta}^{\prime}$ for $P_{\vT}$-almost all ($P_{\vT}$-a.a. for short) $\theta\in\vX$. But the consistency of $\{P_{\theta}\}_{\theta\in\vX}$ together with {(b)} yields that condition (\ref{mag3}) holds true; hence setting $A=X^{-1}(B)$ with $B\in{\mf{B}}$ we deduce 
that $\mathbf{K}(\vT)(B,\cdot)=(P_{\bdot})_X(B)\circ\vT\;$ $P\mid\sigma(\vT)$-a.s..

If no confusion arises, we denote $(P_{\theta})_X$ by $\mathbf{K}(\theta)$ for $\theta\in\vX$.
\end{rems}
{\em Throughout what follows, the conditional distribution $\mathbf{K}(\vT)$ involving in Remark \ref{mag}, {(d)} will 
be considered together with the distributions $\mathbf{K}(\theta)$, for $\theta\in\vX$, associated with $\mathbf{K}(\vT)$ 
as in the above remark, without any additional comments.} 
\smallskip

{\em For the remainder of this section, $\{P_{\theta}\}_{\theta\in\vX}$ is a disintegration of $P$ over $P_{\vT}$ consistent with $\vT$ and $\{X_i\}_{i\in I}$ is a non empty family of $\vS$-$T$-measurable maps from $\vO$ into $\vY$}.
\smallskip

The next result extends Lemma 4.3 from \cite{lm1v}.

\begin{lem}\label{g0}
If $\{k_i\}_{i\in I}$ is a non empty family of $T$-$Z$-Markov kernels, then for 
each $i\in I$ and for any fixed $B\in T$ the following conditions are equivalent:
\begin{enumerate}
\item
$P_{X_i\mid\vT}(B,\cdot)={K}_i(\vT)(B,\cdot)\quad P\mid\sigma(\vT)-\mbox{a.s.}$;
\item
$P_{\theta}(X_i^{-1}(B))=k_i(B,\theta)$ for $P_{\vT}$-a.a. $\theta\in\vX$. 
\end{enumerate}
In particular, the same remains true if $K_i(\vT)(B,\cdot)$ and $k_i(B,\theta)$ are independent of $i$ for all $B\in T$ and $P_{\vT}$-a.a. $\theta\in\vX$.
\end{lem}

{\bf Proof.}
Let us fix on arbitrary $i\in I$. For all $B\in T$ and $D\in Z$ we obtain that
\begin{eqnarray*}
\lefteqn{\int_{\vT^{-1}(D)}P_{X_i\mid\vT}(B,\cdot)dP
=\int_{\vT^{-1}(D)}{K}_i(\vT)(B,\cdot)dP}\\
&\Longleftrightarrow&\int_{\vT^{-1}(D)}\E_P[\chi_{X_i^{-1}(B)}\mid\sigma(\vT)]dP
=\int_{\vT^{-1}(D)}k_i(B,\cdot)\circ{\vT}dP\\
&\stackrel{(\ref{mag3})}{\Longleftrightarrow}&
\int_DP_{\theta}\bigl(X_i^{-1}(B)\bigr)P_{\vT}(d\theta)=\int_Dk_i(B,\theta)P_{\vT}(d\theta).
\end{eqnarray*}
Consequently, the equivalence of assertions $(i)$ and $(ii)$ follows.\hfill$\Box$

\begin{lem}\label{cg0}
Let $\{k_i\}_{i\in I}$ be as in Lemma \ref{g0}. 
Suppose that $I$ is countable and $T$ is countably generated. Then the following are equivalent:
\begin{enumerate}
\item
Condition $P_{X_i\mid\vT}=K_i(\vT)\;$ $P\mid\sigma(\vT)$-a.s. holds true 
for each $i\in I$;
\item
for $P_{\vT}$-a.a. $\theta\in\vX$ condition 
$P_{\theta}\circ X_i^{-1}=k_i(\cdot,\theta)$ 
holds true for each $i\in I$. 
\end{enumerate}
In particular, the same remains true if $K_i(\vT)$ and $k_i(\cdot,\theta)$ are independent of $i$.
\end{lem}

{\bf Proof.} If $(i)$ holds true, we then get by Lemma \ref{g0} that for each $i\in I$ and $B\in T$ condition
$$
P_{\theta}(X_i^{-1}(B))=k_i(B,\theta)\quad\mbox{for $P_{\vT}$-a.a. $\theta\in\vX$}.
$$
is satisfied, which is equivalent to the fact that 
$$
\forall\;i\in{I}\;\;\forall\;B\in{T}\;\;\exists\;{\wt{L}}_{I,i,B}\in Z_0\;\;\forall\;\theta\notin{\wt{L}}_{I,i,B}\quad P_{\theta}(X_i^{-1}(B))=k_i(B,\theta),
$$
where $Z_0:=\{L\in{Z}:P_{\vT}(L)=0\}$. 

Since $I$ is countable and $T$ is countably generated, letting ${\wt{L}}_{I}:=\bigcup_{B\in{\mathcal{G}}_T}\bigcup_{i\in I}{\wt{L}}_{I,i,B}$, where ${\mathcal{G}}_T$ is a countable generator of $T$, and applying a monotone class argument, we find a $P_{\vT}$-null set ${\wt{L}}_{I}\in{Z}$ such that for any $\theta\notin{\wt{L}}_{I}$ the equality $P_{\theta}(X_i^{-1}(B))=k_i(B,\theta)$ holds true.
So assertion $(ii)$ follows.

Applying a similar reasoning we obtain the converse implication.\hfill$\Box$ 
\medskip

The following result extends Lemma 4.1 from \cite{lm1v}. 

\begin{lem}\label{8}
Let $I$ be countable and $T$ countably generated. Then the family $\{X_i\}_{i\in I}$ is
$P$-conditionally independent if and only if for $P_{\vT}$-a.a. $\theta\in\vX$ 
it is $P_{\theta}$-independent.
\end{lem}

{\bf Proof.} Assume that $\{X_i\}_{i\in I}$ is $P$-conditionally independent. 
Then according to Remark \ref{mag}, (b) and following the same reasoning as in 
the proof of \cite{lm1v}, Lemma 4.1, 
we get that
$$
\int_{D}{P_{\theta}}\Bigl(\bigcap_{j=1}^{m}\{X_{i_j}\in B_j\}\Bigr){P_{\vT}}(d\theta)
=\int_{D}\prod_{j=1}^{m}{P_{\theta}}(\{X_{i_j}\in B_j\}){P_{\vT}}(d\theta)
$$
whenever $D\in{Z}$, $m\in\N$, $i_1,\ldots,i_m\in{I}$ are distinct, and $B_1,\ldots,B_m\in{T}$, equivalently that for each $m\in\N$, for all $i_1,\ldots,i_m\in{I}$ distinct and for all $B_1,\ldots,B_m\in{T}$ there exists a $P_{\vT}$-null set ${{L}}_{{I,}m,i_1,\ldots,i_m,B_1,\ldots,B_m}\in{Z}$ such that for any $\theta\notin{{L}}_{{I,}m,i_1,\ldots,i_m,B_1,\ldots,B_m}$ condition 
\begin{equation}\label{ast0}
P_{\theta}\Bigl(\bigcap_{j=1}^{m}\{X_{i_j}\in B_{j}\}\Bigr)
=\prod_{j=1}^{m}{P_{\theta}}(\{X_{i_j}\in B_{j}\})
\end{equation}
holds true. Without loss of generality we may and do assume that $m=2$. Since $T$ is countably generated, applying successively two monotone class arguments we get that there exists a $P_{\vT}$-null set ${L_I}\in{Z}$ such that for any $\theta\notin{L_I}$ condition (\ref{ast0}) holds true for $m=2$, for each $i_1,i_2\in I$ with $i_1\neq i_2$ and for each $B_1,B_2\in{T}$; hence $\{X_i\}_{i\in I}$ is $P_{\theta}$-independent for any $\theta\notin{{L}}_{{I}}$. 
Since the inverse implication is clear, this completes the proof.\hfill$\Box$ 

\begin{lem}\label{id00}
Let $I$ be countable and $T$ countably generated. 
Then the following hold true: 
\begin{enumerate}
\item
The family $\{X_i\}_{i\in I}$ is $P$-conditionally identically distributed if and only if 
for $P_{\vT}$-a.a. $\theta\in\vX$ it is $P_{\theta}$-identically distributed.
\item
The family $\{X_i\}_{i\in I}$ is $P$-conditionally i.i.d. if and only if 
for $P_{\vT}$-a.a. $\theta\in\vX$ it is $P_{\theta}$-i.i.d..
\end{enumerate}
\end{lem}

{\bf Proof.} Ad $(i)$: If $\{X_i\}_{i\in I}$ is $P$-conditionally identically distributed then for any two $i,j\in I$ and for each $B\in{T}$ the equality $P_{X_i\mid\vT}(B)=P_{X_j\mid\vT}(B)$ holds true $P\mid\sigma(\vT)$-a.s..

Then applying successively Remark \ref{mag}, (b) and a monotone class argument as that in the proof of Lemma \ref{cg0}, we find a $P_{\vT}$-null set ${{{\wt{L}}}_{I}^{\prime}}\in{Z}$ such that for any $\theta\notin{{{\wt{L}}}_{I}^{\prime}}$ the family $\{X_i\}_{i\in I}$ is $P_{\theta}$-identically distributed. The inverse implication is immediate by Remark \ref{mag}, (b). 

Ad $(ii)$: Assume that $\{X_i\}_{i\in I}$ is $P$-conditionally i.i.d.. 
It then follows by assertion $(i)$ and Lemma \ref{8} that there exist two $P_{\vT}$-null sets ${{\wt{L}}}_I^{\prime}$ and ${L_I}$ in $Z$ such that for any $\theta\notin{{\wh{L}}_I}:={{\wt{L}}}_I^{\prime}\cup{L_I}$ the family 
$\{X_i\}_{i\in I}$ is $P_{\theta}$-i.i.d.. Since the inverse implication is clear, this completes the proof.\hfill$\Box$

\begin{prop}\label{31}
The counting process $\cnp$ is a $(P,\mathbf{K}(\vT))$-MRP if and only if 
for $P_{\vT}$-a.a. $\theta\in{\R^{d}}$ it is a $(P_{\theta},\mathbf{K}(\theta))$-RP.
\end{prop}

{\bf Proof.} Assume that $\cnp$ is a $(P,\mathbf{K}(\vT))$-MRP, i.e. that the process $\cip$ is $P$-conditionally independent and that for all interarrival times $W_n$ condition $P_{W_n\mid\vT}=\mathbf{K}(\vT)$ holds true $P\mid\sigma(\vT)$-a.s.. Applying now Lemmas \ref{8} and \ref{cg0}, we equivalently get that there exist two $P_{\vT}$-null sets ${H_{\N}}$ and ${{\wt{H}}_{\N}}$ in $Z$ such that for any $\theta\notin{{H}}_{*}:={H_{\N}}\cup{{\wt{H}}_{\N}}$ the sequence $\cip$ is $P_{\theta}$-independent and $(P_{\theta})_{W_n}=\mathbf{K}(\theta)$ for each $n\in\N$, respectively, i.e. such that $\cnp$ is a $(P_{\theta},\mathbf{K}(\theta))$-RP for any $\theta\notin{{H}}_{*}$.\hfill$\Box$

\begin{cor}\label{25p}
Let be given a $\sigma$-subalgebra $\mathcal{F}$ of $\vS$ and a subfield r.c.p. $\{P_{\omega}\}_{\omega\in\vO}$ for $P$ over $R:=P\mid\mathcal{F}$. Then $\cip$ is $P$-conditionally i.i.d. over $\mathcal{F}$ with a conditional probability distribution $\mathbf{K}({\mathrm{id}}_{\vO})=P_{W_n\mid\mathcal{F}}$ $P\mid\mathcal{F}$-a.s. for each $n\in\N$, if and only if for $R$-a.a. $\omega\in\vO$ it is $P_{\omega}$-independent and condition $(P_{\omega})_{W_n}=\mathbf{K}(\omega)$ holds for each $n\in\N$. 
\end{cor}

{\bf Proof.} Put $(\vX,Z):=(\vO,\mathcal{F})$ and $\vT:={\mathrm{id}}_{\vO}$. Then $\{P_{\omega}\}_{\omega\in\vO}$ 
is a disintegration of $P$ over $P_{\vT}=R$ consistent with the map $\vT$. 
So the result follows by Lemma \ref{id00}, $(ii)$.\hfill$\Box$
\smallskip

Finally, we extend Lemmas \ref{8} and \ref{id00} for uncountable index set. 
To this aim, we need to recall some notions more. 

Given a partially ordered set $I$, any increasing family $\{\vS_i\}_{i\in I}$ of
$\sigma$-subalgebras of $\vS$ is said to be a {\bf filtration} for $(\vO,\vS)$. For any family $\{Z_i\}_{i\in I}$ of 
$\vS$-$T$-measurable maps, the filtration $\{{\mathcal Z}_i\}_{i\in I}$ with ${{\mathcal Z}}_i:=\sigma(\bigcup_{j\leq i}\sigma(Z_j))$ for each $i\in I$, is called {\bf the canonical filtration} for $\{Z_i\}_{i\in I}$. 
In particular, for $I=\R_+$ the filtration $\{{\mathcal{Z}}_t\}_{t\in\R_+}$ is said to be 
{\bf right-continuous} if ${\mathcal{Z}}_t=\bigcap_{s>t}{\mathcal{Z}}_s$ for any $t\in\R_+$.

Let $I$ be an arbitrary subset of $\R_+$ and let $\vY$ be a metric space. 
We say that the family $\{X_i\}_{i\in I}$ of $\vS$-$\mf{B}(\vY)$-measurable maps from $\vO$ into $\vY$ is {\bf separable}, if there exists a countable set $G\subseteq I$ such that for each $\omega\in\vO$ the set $\{(u,X_u(\omega)):u\in G\}$ is dense in $\{(i,X_i(\omega)):i\in I\}$. Any such set $G$ is called a {\bf separator} (or {\bf separating set}) for $\{X_i\}_{i\in I}$. 

\begin{rems}\label{phr}
\normalfont
Let $I\subseteq\R_+$ and let $Q$ be a probability measure on $\vS$. 
Then the following can be easily proven:

\noindent{\bf(a)} 
If $\{{U_t}\}_{t\in I}$ is a family of $\vS$-$\mf{B}(\vY)$-measurable maps from $\vO$ into $\vY$, and $\{{{\mathcal{Z}}}_t\}_{t\in I}$ is its canonical filtration, then $\{{U_t}\}_{t\in I}$ is $Q$-independent if and only if for every bounded $\mf{B}(\vY)$-measurable real-valued function $f$ on $\vY$ the equality
\begin{equation}\label{defi}
\E_Q[\chi_Af({U_t})]=Q(A)\E_Q[f({U_t})]
\end{equation}
holds true for each $s,t\in I$ with $s<t$ and for each $A\in{{\mathcal{Z}}}_s$. 

\noindent{\bf(b)}
If $U_1$ and $U_2$ are two $\vS$-$\mf{B}(\vY)$-measurable maps from $\vO$ into $\vY$, then they are $Q$-identically distributed if and only if $\E_Q[f(U_1)]=\E_Q[f(U_2)]$ 
for every bounded $\mf{B}(\vY)$-measurable real-valued function $f$ on $\vY$.
\end{rems}
Recall that the family $\{X_t\}_{t\in\R_+}$ has {\bf $P$-(conditionally) independent increments}, if for each $m\in\N$ and for each $t_0,t_1,\ldots,t_m\in\R_+$, such that $0=t_0<t_1<\cdots<t_m$ the increments $X_{t_j}-X_{t_{j-1}}$ ($j\in\N_m$) are $P$-(conditionally) independent.

\begin{prop}\label{nid0}
Let $\vY$ be a Polish space, let $\{X_t\}_{t\in\R_+}$ be a family of
$\vS$-$\mf{B}(\vY)$-measurable maps from $\vO$ into $\vY$ and let $\{\mathcal{H}_t\}_{t\in\R_+}$ be its canonical filtration. 
If the family $\{X_t\}_{t\in\R_+}$ is separable with separator $\Q_+$, then the following hold true:
\begin{enumerate}
\item
If $\{{\mathcal{H}}_t\}_{t\in\R_+}$ is right-continuous, then $\{X_t\}_{t\in\R_+}$ is $P$-conditionally independent if and only if for $P_{\vT}$-a.a. $\theta\in\vX$ 
it is $P_{\theta}$-independent.
\item
The family $\{X_t\}_{t\in\R_+}$ is $P$-conditionally identically distributed if and only if for $P_{\vT}$-a.a. $\theta\in\vX$ it is $P_{\theta}$-identically distributed.
\item
If $\{{\mathcal{H}}_t\}_{t\in\R_+}$ is right-continuous, then $\{X_t\}_{t\in\R_+}$ is $P$-conditionally i.i.d. if and only if for $P_{\vT}$-a.a. $\theta\in\vX$ it is $P_{\theta}$-i.i.d..
\item
If $\vY$ is either an open or a closed subset of $\ov\R$ and if $\{{\wt{\mathbf{K}}_t}\}_{t\in\R_+}$ is a family of $\mf{B}(\vY)$-$Z$-Markov kernels such that ${\wt{\mathbf{K}}_t}(\theta)$ is for every $t\in\R_+$ a probability distribution on $\mf{B}(\vY)$ with parameter $\theta\in\vX$, and the function $t\longmapsto{\wt{\mathbf{K}}_t}(\theta)(B)$ is continuous for any fixed $B\in\mf{B}(\vY)$ and $\theta\in\vX$, then condition $P_{X_t\mid\vT}={\wt{\mathbf{K}}_t}(\vT)$  $P\mid\sigma(\vT)$-a.s. holds true for each $t\in\R_+$ if and only if for $P_{\vT}$-a.a. $\theta\in\vX$ condition $P_{\theta}\circ{X_t}^{-1}=\wt{\mathbf{K}}_t(\theta)$ holds true for each $t\in\R_+$.
\end{enumerate}
Moreover, assertions $(i)$ to $(iii)$ remain true for the increments of $\{X_t\}_{t\in\R_+}$ in the place of $\{X_t\}_{t\in\R_+}$.
\end{prop}

{\bf Proof.} Ad $(i)$: The ``if'' implication follows as in the proof of Lemma 4.1 from \cite{lm1v}. For the ``only if'' part, assume that $\{X_t\}_{t\in\R_+}$ is $P$-conditionally independent and note that our assumptions for $\{X_t\}_{t\in\R_+}$ and $\{\mathcal{H}_t\}_{t\in\R_+}$ imply 
${\mathcal{H}}_s=\sigma(\{X_u\}_{u\in\Q_+,u\leq s})
=\bigcap_{s^{\prime}\in\Q_+,s^{\prime}>s}\mathcal{H}_{s^{\prime}}$ for each $s\in\R_+$.

{\bf(a)} It follows by Lemma \ref{8} that there exists a $P_{\vT}$-null set $O_{\Q_+}\in{Z}$ such that for any $\theta\notin{O}_{\Q_+}$ condition (\ref{ast0}) holds true with $\Q_+$ and $\mf{B}(\vY)$ in the place of $I$ and $T$, respectively. 

Throughout this proof fix on an arbitrary $\theta\notin{O}_{\Q_+}$. Then condition (\ref{ast0}) together with Remark \ref{phr}, (a) implies that for all $s,t\in\Q_+$ with $s<t$, for every bounded $\mf{B}(\vY)$-measurable real-valued function $f$ on $\vY$ and for each $A\in\mathcal{H}_s$ we have
\begin{equation}\label{dalpha}
\E_{P_{\theta}}[\chi_Af(X_t)]=P_{\theta}(A)\E_{P_{\theta}}[f(X_t)].
\end{equation}
If we take $s,t\in\R_+$ with $s<t$ and if we write (\ref{dalpha}) 
for $s^{\prime},t^{\prime}\in\Q_+$ with $s^{\prime}<t^{\prime}$ and then let $s^{\prime}\downarrow s$ and $t^{\prime}\downarrow t$, the separability of $\{X_t\}_{t\in\R_+}$ together with an application of Lebesgue's Dominated Convergence Theorem yields that for all $A\in\bigcap_{s^{\prime}\in\Q_+,s^{\prime}>s}\mathcal{H}_{s^{\prime}}=\mathcal{H}_s$ 
and for every bounded continuous real-valued function $f$ on $\vY$ condition (\ref{dalpha}) holds true.

{\bf(b)} 
Let $s,t\in\R_+$ with $s<t$ and let $f$ be a function as in (\ref{defi}). Then for each $n\in\N$ there exists a bounded continuous real-valued function $g_n$ on $\vY$ satisfying the inequality $\int|g_n-f|d(P_{\theta})_{X_t}\leq\frac{1}{n}$ (cf. e.g. \cite{fr4}, Proposition 415P); hence there exists a sequence $\{g_n\}_{n\in\N}$ of bounded continuous real-valued functions on $\vY$ such that condition $\lim_{n\to\infty}\int\chi_A(|g_n-f|\circ X_t)dP_{\theta}=0$ 
holds true for all $A\in\mathcal{F}_s$, implying together with (a) that 
$\E_{P_{\theta}}[\chi_Af(X_t)]=\lim_{n\to\infty}\E_{P_{\theta}}[\chi_Ag_n(X_t)]
=\lim_{n\to\infty}P_{\theta}(A)\E_{P_{\theta}}[g_n(X_t)]
=P_{\theta}(A)\E_{P_{\theta}}[f(X_t)]$ for all $A\in\mathcal{F}_s$; hence by Remark \ref{phr}, (a) we get that $\{X_t\}_{t\in\R_+}$ is $P_{\theta}$-independent, which proves $(i)$.

Ad $(ii)$: The ``if'' implication is immediate by Remark \ref{mag}, (b). For the ``only if'' part, assume that $\{X_t\}_{t\in\R_+}$ is $P$-conditionally identically distributed. 

\noindent{\bf{(c)}}
Since $\{X_t\}_{t\in\R_+}$ is $P$-conditionally identically distributed, we get that 
for any two $s,t\in\Q_+$ and for each $B\in\mf{B}(\vY)$ the equality $P_{X_t\mid\vT}(B)=P_{X_s\mid\vT}(B)$ holds $P\mid\sigma(\vT)$-a.s.. The latter together with Lemma \ref{id00}, $(i)$ yields the existence of a $P_{\vT}$-null set 
${\wt{O}}_{\Q_+}\in{Z}$ such that for any $\theta\notin{\wt{O}}_{\Q_+}$ and for all $s,t\in\Q_+$ condition $(P_{\theta})_{X_t}=(P_{\theta})_{X_s}$ holds true, which by Remark \ref{phr}, (b) equivalently yields that for any $\theta\notin{\wt{O}}_{\Q_+}$, for every function $f$ as in the above remark, and for all $s,t\in\Q_+$ we have $\E_{P_{\theta}}[f(X_t)]=\E_{P_{\theta}}[f(X_s)]$.

Till the end of the proof of $(ii)$, fix on an arbitrary $\theta\notin{\wt{O}}_{\Q_+}$.

\noindent{\bf{(d)}}
If we take $s,t\in\R_+$ and if we write the last equality
for $s^{\prime},t^{\prime}\in\Q_+$ and then let $s^{\prime}\downarrow s$ and $t^{\prime}\downarrow t$, the separability of $\{X_t\}_{t\in\R_+}$ together with an application of Lebesgue's Dominated Convergence Theorem yields that for 
every bounded continuous real-valued function $f$ on $\vY$ condition $\E_{P_{\theta}}[f(X_t)]=\E_{P_{\theta}}[f(X_s)]$ holds true.

Following now the same reasoning with that of steps (b) and (c), 
we obtain that the last equality is satisfied by all functions $f$ as in Remark \ref{phr}, and all $s,t\in\R_+$, which is equivalent to the fact that condition $(P_{\theta})_{X_t}=(P_{\theta})_{X_s}$ holds true for all $s,t\in\R_+$; hence $(ii)$ follows.

Ad $(iii)$: Assume that $\{X_t\}_{t\in\R_+}$ is $P$-conditionally i.i.d. and that its canonical filtration is right-continuous. It then follows by assertions $(i)$ and $(ii)$ there exist two $P_{\vT}$-null sets ${\wt{O}}_{\Q_+}$ and $O_{\Q_+}$ in $Z$ such that for any $\theta\notin{\wh{O}}_{\Q_+}:={\wt{O}}_{\Q_+}\cup{O}_{\Q_+}$ the family 
$\{X_t\}_{t\in\R_+}$ is $P_{\theta}$-i.i.d.. Since the inverse implication is clear, 
assertion $(iii)$ follows.

Ad $(iv)$: Assume that condition $P_{X_t\mid\vT}={\wt{\mathbf{K}}_t}(\vT)$ holds $P\mid\sigma(\vT)$-a.s. for any $t\in\R_+$. It then follows by Lemma \ref{cg0} that there exists a $P_{\vT}$-null set $\wt{O}_{\Q_+}^{\prime}\in{Z}$ such that for any $\theta\notin{\wt{O}}_{\Q_+}^{\prime}$, for each $B\in\mf{B}(\vY)$ and for any $t\in\Q_+$ the following condition holds true: 
\begin{equation}\label{axr}
P_{\theta}(X_t^{-1}(B))={\wt{\mathbf{K}}_t}(\theta)(B). 
\end{equation}
Fix on arbitrary $\theta\notin{{\wt{O}}_{\Q_+}^{\prime}}$ and $t\in\R_+$. Then the separability of $\{X_t\}_{t\in\R_+}$ implies that there exists a monotone sequence $\{X_s\}_{s\in\Q_+}$ such that $s\to{t}$ and $X_t=\lim_{s\to{t}}X_s$, which together with (\ref{axr}) and the Monotone Convergrence Theorem yields that $$(P_{\theta})_{X_t}=\lim_{s\to{t}}(P_{\theta})_{X_s}
=\lim_{s\to{t}}\wt{\mathbf{K}}_s(\theta)=\wt{\mathbf{K}}_{t}(\theta).$$ 
Since the inverse implication follows by applying similar arguments, we obtain $(iv)$.
\smallskip

Moreover, the proofs of assertions $(i)$ to $(iii)$ for the increments of $\{X_t\}_{t\in\R_+}$ run in the same way as for $\{X_t\}_{t\in\R_+}$.\hfill$\Box$
\medskip

It is immediate from the corresponding definitions that if $\{X_t\}_{t\in\R_+}$ satisfies condition $X_0(\omega)=0$ for each $\omega\in\vO$ and has $P$-conditionally independent increments, then it will have $P$-conditionally stationary increments, if and only if for each $t,h\in\R_+$ the equality $P_{X_{t+h}-X_t\mid\vT}=P_{X_h\mid\vT}$ holds $P\mid\sigma(\vT)$-a.s. true (cf. e.g. \cite{lm1v}, page 68 for the definition of conditionally stationary increments).
\smallskip

The following result extends a basic one from \cite{lm1v}, that is Proposition 4.4. 

\begin{thm}\label{19}
Let $\cnp$ be a counting process and let $\{{\wt{\mathbf{K}}_t}\}_{t\in\R_+}$ 
be as in Proposition \ref{nid0} but with $\vY=[0,\infty]$. 
Then $\cnp$ has $P$-conditionally stationary independent increments such that condition 
$$
P_{N_t\mid\vT}={\wt{\mathbf{K}}_t}(\vT)\quad{P}\mid\sigma(\vT)-\mbox{a.s.}
$$ 
holds true for each $t\in\R_+$ if and only if for $P_{\vT}$-a.a. $\theta\in\vX$ it 
has $P_{\theta}$-stationary independent increments such that $(P_{\theta})_{N_t}={\wt{\mathbf{K}}_t}(\theta)$ for each $t\in\R_+$.
\end{thm}

{\bf Proof.} Since $\cnp$ is a counting process it has right-continuous paths; hence it is separable with separator $\Q_+$. Note also that the canonical filtration of $\cnp$ is right-continuous (see \cite{pro}, Theorem 25, where the proof works for any probability space not necessarily complete). Thus, all assumptions of Proposition \ref{nid0} are fulfilled, and so we may apply it to deduce the thesis of the theorem.\hfill$\Box$

\begin{cor}[\cite{lm1v}, Proposition 4.4]\label{c19}
The family $\cnp$ is a $P$-MPP with parameter $\vT$ if and only if it is a $P_{\theta}$-Poisson process with parameter $\theta$ for $P_{\vT}$-a.a. $\theta{\in\R}$.
\end{cor}

\section{Mixed renewal processes and exchangeability}\label{exch}

An infinite family $\{X_i\}_{i\in I}$ of $\vS$-$T$-measurable maps from  $\vO$ into $\vY$ is said to be {\bf exchangeable} under $P$ or {\bf $P$-exchangeable}, if for each $r\in\N$ we have 
$$
P\Bigl(\bigcap_{k=1}^{r}X_{i_k}^{-1}(B_k)\Bigr)
=P\Bigl(\bigcap_{k=1}^{r}X_{j_k}^{-1}(B_k)\Bigr)
$$
whenever $i_1,\ldots,i_r\in{I}$ are distinct, $j_1,\ldots,j_r\in{I}$ are distinct, and $B_k\in T$ for each $k\leq{r}$ (cf. e.g. \cite{fr4}, 459C).

For the purposes of this section we recall the notions of product r.c.p. as well as of infinite products of measure spaces.

Let $Q$ be a probability measure on $T$. Assume that $M$ is a probability on the $\sigma$-algebra $\vS\otimes{T}$ such that $P$ and $Q$ are the marginals of $M$. Assume also that for each $y\in\vY$ there exists a probability $P_{y}$ on $\vS$, satisfying the following properties:
\begin{description}
\item[(D1)] For every $A\in\vS$ the map $y\longmapsto P_{y}(A)$ is $T$-measurable;
\item[(D2)] $M(A\times B)=\int_BP_{y}(A)Q(dy)$ for each $A\times B\in\vS\times{T}$.
\end{description}
Then, $\{P_{y}\}_{y\in\vY}$ is said to be a {\bf product r.c.p.} {\rm on $\vS$ for $M$ with respect to $Q$} (see \cite{fa}, Section 2 or \cite{mms3}, Definition 1.1).
\smallskip
 
Let $I$ be an arbitrary non-empty index set. 
If $\{(\vO_{i},\vS_{i},P_{i})\}_{i\in I}$ is a family of probability spaces then, for each $\emptyset\neq J\subseteq I$ we denote by $(\vO_{J},\vS_{J},P_{J})$ the product probability space $\otimes_{i\in J}(\vO_i,\vS_i,P_i)
:=(\prod_{i\in J}\vO_i,\otimes_{i\in J}\vS_i,\otimes_{i\in J}P_i)$. 
If $(\vO,\vS,P)$ is a probability space, we write $P_I$ for the product measure on $\vO^I$ and $\vS_I$ for its domain.

\begin{lem}\label{20}
Let $\mathcal{F}$ be a $\sigma$-subalgebra of $\vS$ and let $\{X_i\}_{i\in I}$ be a non empty family of $\vS$-$T$-measurable maps from $\vO$ into $\vY$ such that $\{X_i\}_{i\in I}$ is $P$-conditionally i.i.d. over $\mathcal{F}$.
Suppose that $T$ is countably generated and that $P_{X_i}$ is perfect for each $i\in I$. 
Then there exists a probability measure $M$ on $T\otimes\mathcal{F}$ with marginal $R:=P\mid\mathcal{F}$ on $\mathcal{F}$ such that $M:=P\circ(X_i\times{\mathrm{id}}_{\vO})^{-1}$ for every $i\in I$, and a product r.c.p. $\{Q_{\omega}\}_{\omega\in\vO}$ on $T$ for $M$ with respect to $R$, such that
\begin{enumerate}
\item
for any fixed $B\in T$ and $i\in{I}$ the map $Q_{\bdot}(B):\vO\longrightarrow[0,1]$ is $R$-a.s. equal to $P(X_{i}^{-1}(B)\mid\mathcal{F})(\cdot)$;
\item
$\int_{F}Q_{\omega}^{I}(H)R(d\omega)=P(F\cap{X^{-1}}(H))$ for every $F\in\mathcal{F}$ and 
$H\in T_{I}$, where $Q_{\omega}^{I}$ denotes the $I$-fold product probability 
$\otimes_{i\in I}P_i$ of copies $P_i:=Q_{\omega}$ of $Q_{\omega}$ for 
$i\in I$, and $X:\vO\longrightarrow\vY^{I}$ is defined by $X(\omega)=\bigl(X_i(\omega)\bigr)_{i\in I}$ for each $\omega\in\vO$.
\end{enumerate}
\end{lem}

{\bf Proof.} First fix on an arbitrary $i\in I$.

\noindent{\bf(a)} 
The function 
$X_i\times{\mathrm{id}}_{\vO}$ from $\vO$ into $\vY\times\vO$ defined by means of 
$$
(X_i\times{\mathrm{id}}_{\vO})(\omega):=(X_i(\omega),\omega)\quad\mbox{for each}\;\;\omega\in\vO
$$ 
is $\vS$-${T}\otimes{\mathcal{F}}$-measurable. So, we have a probability measure 
$M_i:=P\circ(X_i\times{\mathrm{id}}_{\vO})^{-1}$ on ${T}\otimes{\mathcal{F}}$. Since 
all $X_i$ are $P$-conditionally identically distributed over $\mathcal{F}$, it 
follows that $M_i$ is independent of $i$, so we may write $M:=M_{i^*}$ for any fixed $i^{*}\in I$.
\smallskip

\noindent{\bf(b)} 
There exists a product r.c.p. 
$\{Q_{\omega}\}_{\omega\in\vO}$ on $T$ for $M$ with respect to $R=P\mid\mathcal{F}$ 
such that for any fixed $B\in T$
$$
Q_{\bdot}(B)=P(X_i^{-1}(B)\mid\mathcal{F})(\cdot)\quad{R}-\mbox{a.s.}. 
$$
In fact, by assumption each marginal measure $P_{X_i}$ of $M$ on $T$ is perfect and $T$ is countably generated; hence by \cite{fa}, Theorem 6, there exists a product r.c.p. $\{Q_{\omega}\}_{\omega\in\vO}$ on $T$ for $M$ with respect to $R$. 

Since $\{Q_{\omega}\}_{\omega\in\vO}$ satisfies (D2), we get that
$$
\int_{F}Q_{\omega}(B)R(d\omega)=M(B\times F)=P(F\cap X_i^{-1}(B))
=\int_{F}P(X_i^{-1}(B)\mid\mathcal{F})(\omega)R(d\omega)
$$
for every $B\in{T}$ and $F\in\mathcal{F}$, which proves (b); hence $(i)$ follows. 

\noindent{\bf(c)} 
Using $(i)$ and a monotone class argument we get that $(ii)$ holds true.\hfill$\Box$ 
\medskip

The next result extends a corresponding one due to Olshen (see \cite{ols}, Theorem (3)) concerning a generalization of de Finetti's Theorem.

\begin{prop}\label{bolg}
Let $\{X_i\}_{i\in I}$ be a $P$-exchangeable infinite family of $\vS$-$T$-measurable maps from $\vO$ into $\vY$. Suppose that $T$ is countably generated and $P_{X_i}$ is perfect for each $i\in{I}$. Then there exists a $d$-dimensional random vector $\vT$ such that $\{X_i\}_{i\in{I}}$ is $P$-conditionally i.i.d. given $\vT$.
\end{prop}

{\bf Proof.}  
Since $\{X_i\}_{i\in I}$ is $P$-exchangeable, it follows by \cite{fr4}, Theorem 459B, that there exist a $\sigma$-subalgebra $\mathcal{F}$ of $\vS$ such that $\{X_i\}_{i\in{I}}$ is $P$-conditionally i.i.d. over $\mathcal{F}$. So, applying Lemma \ref{20}, we deduce that there exists a family $\{Q_{\omega}\}_{\omega\in\vO}$ of $T$-$\mathcal{F}$-Markov kernels such that
\begin{equation}\label{bo1}
\int_{F}Q_{\omega}^{I}(H)R(d\omega)=P(F\cap X^{-1}(H))
\end{equation}
for every $H\in T_{I}$ and $F\in\mathcal{F}$, where $R:=P\mid\mathcal{F}$. 
Then there exists a countably generated $\sigma$-subalgebra $\mathcal{A}$ of $\mathcal{F}$ such that $Q_{\bullet}(B)$ is $\mathcal{A}$-measurable for arbitrary but fixed $B\in{T}$ (take e.g. 
$\mathcal{A}_{{B}}:=\{[Q_{\bullet}(B)]^{-1}(E): E\in{\mathcal{G}}_{\mf{B}}\}$ for $B\in{T}$, and 
$\mathcal{A}:=\sigma(\bigcup_{B\in{\mathcal{G}}_T}\mathcal{A}_B)$, where $\mathcal{G}_{\mf{B}}$ and $\mathcal{G}_{T}$ is a countable generator of $\mf{B}$ and $T$, respectively). Since $\mathcal{A}$ is countably generated, there exists a map $\wt{\vT}:\vO\longrightarrow\R$ such that $\mathcal{A}=\sigma(\wt{\vT})$ (take e.g. $\wt{\vT}$ to be the {\em Marczewski functional} on $\vO$, cf. e.g. \cite{fr4}, 343E for the definition). But since $\{X_i\}_{i\in I}$ is $P$-conditionally i.i.d. over $\mathcal{F}$ and $\mathcal{A}\subseteq\mathcal{F}$, it follows that $\{X_i\}_{i\in I}$ is so over $\mathcal{A}=\sigma(\wt{\vT})$.

Note also that $\R$ and $\R^{d}$ are standard Borel spaces of the same cardinality, there exists a Borel isomorphism $g$ from $\R$ into $\R^d$ (cf. e.g. \cite{fr4}, Corollary 424D(a)). So, putting $\vT:=g\circ\wt{\vT}$, we get that $\vT$ is a $d$-dimensional random vector on $\vO$ such that $\sigma(\vT)=\sigma(\wt{\vT})$.\hfill$\Box$

\begin{cor}[see Olshen, R. \cite{ols}, Theorem (3)]\label{bcor}
If $\{X_n\}_{\in\N}$ is a $P$ - exchangeable sequence of measurable maps from $\vO$ into 
a complete, separable metric space, then there exists a real-valued random variable $\vT$ 
on $\vO$ such that $\{X_n\}_{n\in\N}$ is $P$-conditionally i.i.d. given $\vT$.
\end{cor}

\begin{thm}\label{21}
Let $\{X_i\}_{i\in I}$ be an infinite family of $\vS$-$T$-measurable maps from $\vO$ into $\vY$. Consider the following assertions:
\begin{enumerate}
\item 
$\{X_i\}_{i\in I}$ is $P$-exchangeable.
\item
There exists a $\sigma$-subalgebra $\mathcal{F}$ of $\vS$ such that 
$\{X_i\}_{i\in I}$ is $P$-conditionally i.i.d. over $\mathcal{F}$.
\item
There exists a $\sigma$-subalgebra $\mathcal{F}$ of $\vS$ and a family $\{Q_{\omega}\}_{\omega\in\vO}$ of $T$-$\mathcal{F}$-Markov kernels such that
$$
\int_{F}Q_{\omega}^{I}(H)R(d\omega)=P(F\cap X^{-1}(H))
$$
for every $H\in T_{I}$ and $F\in\mathcal{F}$, where $R:=P\mid\mathcal{F}$ and $Q_{\omega}^{I}$, $X$ are as in Lemma \ref{20}.
\item
There exists a $\vS$-$\mf{B}_d$-measurable map $\vT$ from $\vO$ into $\R^d$ such that $\{X_i\}_{i\in{I}}$ is $P$-conditionally i.i.d. given $\vT$.
\end{enumerate}
Then $(i)\Longleftrightarrow (ii)$, $(iii)\Longrightarrow(i)$ and $(iv)\Longrightarrow(i)$. 
If any one of conditions $(i)$ to $(iv)$ is satisfied, then all image measures $P_{X_i}$ are equal.

Moreover, if $P_{X_i}$ is perfect for any $i\in I$ and $T$ is countably generated, then 
assertions $(i)$ to $(iv)$ are equivalent. 
\end{thm}

{\bf Proof.} 
The equivalence $(i)\Longleftrightarrow(ii)$ follows by \cite{fr4}, Theorem 459B, while the implications $(iii)\Longrightarrow(i)$ and $(iv)\Longrightarrow(i)$ are evident.

Clearly, if assertion $(i)$ or equivalently $(ii)$ is satisfied then all $P_{X_i}$ are equal and the same applies if $(iii)$ or $(iv)$ holds true.
\smallskip

Moreover, if every measure $P_{X_i}$ is perfect and $T$ is countably generated, then implications $(ii)\Longrightarrow(iii)$ and $(i)\Longrightarrow(iv)$ follow from Lemma \ref{20} and Proposition \ref{bolg}, respectively. So we get that assertions $(i)$ to $(iv)$ are equivalent.\hfill$\Box$

\begin{cor}\label{21a}
Let $\{X_t\}_{t\in\R_+}$ be a family of $\vS$-$T$-measurable maps from $\vO$ into $\vY$. Suppose that $\vS$ is countably generated, $P$ is perfect, $\vY$ is a Polish space, $\{X_t\}_{t\in\R_+}$ is separable with separator $\Q_+$ and that its canonical filtration is right-continuous. Then each of the items $(i)$ to $(iv)$ of Theorem \ref{21} is equivalent to condition 
\begin{enumerate}\addtocounter{enumi}{4}
\item there exist a $d$-dimensional random vector $\vT$ and a disintegration $\{P_{\theta}\}_{\theta\in\R^d}$ of $P$ over $P_{\vT}$ consistent with $\vT$ such that $\{X_t\}_{t\in\R_+}$ is $P_{\theta}$-i.i.d. for $P_{\vT}$-a.a. $\theta\in\R^d$.
\end{enumerate}
\end{cor}

{\bf Proof.} It follows by Remark \ref{mag}, (a), that given a $d$-dimensional random vector $\vT$, there exists a disintegration $\{P_{\theta}\}_{\theta\in\R^d}$ of $P$ over $P_{\vT}$ consistent with $\vT$. Thus we may apply Proposition \ref{nid0} to obtain that condition $(v)$ is equivalent to $(iv)$ of Theorem \ref{21}. The equivalence of all items $(i)$ to $(v)$ is immediate by Theorem \ref{21}.\hfill$\Box$

\begin{rems}\label{r21}
\normalfont
\noindent{\bf (a)}
The assumption ``$P_{X_i}$ perfect'' made in the last theorem is easily verified in the usual applications, since this is covered by the following facts: 
($\alpha$) If $\vY$ is a Polish space then each $P_{X_i}$ is Radon (cf. e.g. \cite{fr4}, Proposition 434K(b) and \cite{fr4}, Definition 411H(b) for the definition of a Radon measure); hence perfect (cf. e.g. \cite{fr4}, Proposition 416W(a)). 
($\beta$) If $P$ is perfect then each $P_{X_i}$ is so 
(cf. e.g. \cite{fr4}, Proposition 451E(a)). 
($\gamma$) If $\vO=\vY^I$, $X_i\;(i\in{I})$ are the canonical projections from $\vO$ onto $\vY$, and $P$ is any probability measure on $\vS:=T_I$ then each $P_{X_i}$ is perfect if and only if $P$ is perfect (cf. e.g. \cite{fr4}, Theorem 454A(b)(iii)).

\noindent{\bf (b)} 
To the best of our knowledge, the most general result concerning the equivalence of assertions $(i)$ and $(iii)$ of Theorem \ref{21} is Theorem 1.1 from \cite{ka} (which extends de Finetti's Theorem), saying that for each infinite sequence $\{X_n\}_{n\in\N}$ of random variables taking values in a standard Borel space $\vY$ (i.e. $\vY$ is isomorphic to some Borel-measurable subset of $\R$) assertions $(i)$ and $(iii)$ of Theorem \ref{21} with $\{X_n\}_{n\in\N}$ in the place of $\{X_i\}_{i\in{I}}$ are equivalent. It is well-known that any Polish space is standard Borel; in particular, $\R^{d}$ and $\R^{\N}$ are such spaces. 

\noindent{\bf (c)}
There are measurable spaces $(\vY,T)$ satisfying the assumptions of Theorem \ref{21}, i.e. that $T$ is countably generated and each $P_{X_i}$ is perfect, which are not standard Borel spaces; hence Theorem \ref{21} extends Theorem 1.1 from \cite{ka}. In fact, it is known that each uncountable analytic Hausdorff space (i.e. a non-empty topological Hausdorff space being a continuous image of the space $\N^{\N}$, cf. e.g. \cite{fr4}, Definition 423A) has a non-Borel analytic subset (cf. e.g. \cite{fr4}, Proposition 423L). It is also known that for each analytic Hausdorff space $\vY$ the Borel $\sigma$-algebra $\mf{B}(\vY)$ is countably generated (cf. e.g. \cite{fr4}, 423X(d)), and that any Borel probability measure on $\mf{B}(\vY)$ is always inner regular with respect to compact sets (see \cite{hj}, Chapter IV, Theorem 1, page 195); hence it is perfect (cf. e.g. \cite{fr4}, Proposition 451C). Consequently, each uncountable analytic Hausdorff space has a subset satisfying the assumptions of Theorem \ref{21}, but not being a standard Borel space.

\noindent{\bf (d)}
Restricting attention to measurable spaces $(\vY,T)$ satisfying the countability assumption of $T$ as in Theorem \ref{21}, costs us some generality; for instance, the general compact Hausdorff space does not satisfy the countability assumption for $T$, and it is known that 
the equivalence of assertions $(i)$ to $(iii)$ of Theorem \ref{21} is true for countable products of compact Hausdorff spaces (see \cite{hs} or \cite{df}). More general, the equivalence of assertions $(i)$ to $(iii)$ of Theorem \ref{21} is proven in \cite{fr4}, Theorem 459G for uncountable products of general Hausdorff spaces.
But all the above are specialized in the product situation of topological spaces, while assertions $(i)$ to $(iii)$ of Theorem \ref{21} have the advantage of being free from any topological assumption as well as from any product situation.
\end{rems}
The following definition of an MRP traces back to Huang \cite{hu}, Section 1, Definition 3.

\begin{df}\label{hua}
\normalfont
The counting process $\cnp$ is said to be a {\bf $\nu$-mixed renewal process} associated with $\{P_{\wt{y}}\}_{\wt{y}\in\wt\vY}$, if for every $r\in\N$ and for 
all $w_1,\ldots,w_r\in\R$ condition
$$
P\Bigl(\bigcap_{k=1}^{r}\{W_k\leq w_k\}\Bigr)
=\int\prod_{k=1}^{r}P_{\wt{y}}(W_k\leq w_k)\nu(d\wt{y}), 
$$
holds true, where $\{P_{\wt{y}}\}_{\wt{y}\in\wt\vY}$ is a family of probability measures on $\vS$ and $\nu$ is a probability measure on $B(\wt\vY):=\sigma(\{P_{\bdot}(E):E\in\vS\})$ such that for $\nu$-a.a. $\wt{y}\in{\wt{\vY}}$ the process $\{W_n\}_{n\in\N}$ is $P_{\wt{y}}$-identically distributed.
\end{df}
In Huang's definition it is assumed that $\cnp$ takes values only in $\N_0$, which is equivalent to the mild assumption that $\cnp$ has zero probability of explosion, that is 
$P(\bigcup_{t\in(0,\infty)}\{N_t=\infty\})=0$ (cf. e.g. \cite{sch}, Lemma 2.1.4).

\begin{rem}\label{27}
\normalfont
Note that in Huang's \cite{hu} definition the assumption that $\cip$ is $P_{\wt{y}}$-identically distributed for $\nu$-a.a. $\wt{y}\in{\wt{\vY}}$ is not written explicitly. But this assumption must be included there, since it is necessary for the validity of the Corollary on page 20 of \cite{hu}, as it follows from Example \ref{16} below.

In fact, consider the process $\cnp$ of Example \ref{16}, where the above assumption does not hold true, as well as Huang's definition of a $\nu$-MRP without the above assumption. Also note that $q:=P(Z<\infty)=0<1$, where $Z$ is the almost sure limit of $\cnp$ as $t\to\infty$. Assume, if possible, that Corollary in \cite{hu}, page 20, holds true. Then conditional on the event $\{Z=\infty\}$ the process $\cnp$ has the exchangeable property $(E)$ (see \cite{hu}, Definition 1 for the definition) implying that $\cip$ is exchangeable, 
a contradiction to Example \ref{16}.
\end{rem}

\begin{thm}\label{9}
Consider the following assertions:
\begin{enumerate}
\item
There exists a $d$-dimensional ($d\in\N$) random vector $\vT$ such that $\cnp$ is a $(P,\mathbf{K}(\vT))$-MRP.
\item
There exist a random vector $\vT$, a disintegration $\{P_{\theta}\}_{\theta\in{\R^{d}}}$ of $P$ over $P_{\vT}$
consistent with $\vT$, and a family $\{\mathbf{K}(\theta)\}_{\theta\in{\R^{d}}}$ of $\mf{B}$-$Z$-Markov kernels 
for $P_{\vT}$-a.a. $\theta\in\R^d$ the family $\cnp$ is a $(P_{\theta},\mathbf{K}(\theta))$-RP.
\item 
The process $\cip$ is $P$-exchangeable.
\item
There exist a $\sigma$-subalgebra $\mathcal{F}$ of $\vS$ and a family
$\{Q_{\omega}\}_{\omega\in\vO}$ of $\mf{B}$-$\mathcal{F}$-Markov kernels such that
$$
\int_{F}Q_{\omega}^{\N}(H)R(d\omega)=P(F\cap W^{-1}(H))
$$
for every $H\in{\mf{B}}_{\N}$ and $F\in\mathcal{F}$, where $R:=P\mid\mathcal{F}$, $W:=(W_1,\ldots,W_n,\ldots)$ and 
$Q_{\omega}^{\N}$ denotes the $\N$-fold product probability $\otimes_{n\in\N}P_n$ of copies $P_n:=Q_{\omega}$ of $Q_{\omega}$ for $n\in\N$.
\item
There exist a $\sigma$-subalgebra $\mathcal{F}$ of $\vS$, a subfield r.c.p. $\{S_{\omega}\}_{\omega\in\vO}$ for $P$ over the restriction $R:=P\mid\mathcal{F}$, and a family $\{\mathbf{K}(\omega)\}_{\omega\in\vO}$ of $\mf{B}$-$\mathcal{F}$-Markov kernels such that for $R$-a.a. $\omega\in\vO$ the sequence $\cip$ is $S_{\omega}$-independent and condition $(S_{\omega})_{W_n}=\mathbf{K}(\omega)$ holds for each $n\in\N$.
\item
There exist a set ${\wt{\vY}}$, a family $\{S_{\wt{y}}\}_{\wt{y}\in\wt\vY}$ of probability measures on $\vS$ and a probability measure $\nu$ on $B({\wt{\vY}}):=\sigma(\{S_{\bdot}(E): E\in\vS\})$ such that 
$\cnp$ is a $\nu$-MRP associated with $\{S_{\wt{y}}\}_{\wt{y}\in\vY}$.
\end{enumerate}
Then the following implications hold true:
\begin{center} 
\begin{tabular}[t]{ccccc}
$(i)$ & $\Longleftarrow$ & $(ii)$ & $\Longrightarrow$ & $(v)$ 
\\
$\Updownarrow$ & $\;$  & $\Downarrow$ & $\;$ & $\Downarrow$ 
\\
$(iii)$ & $\Longleftrightarrow$ & $(iv)$& $\Longleftarrow$ & $(vi)$ 
\end{tabular}
\end{center}
Moreover, if $\vS$ is countably generated and $P$ is perfect then items $(i)$ to $(vi)$ are all equivalent.
\end{thm}

{\bf Proof.} First note that the implications $(i)\Longrightarrow(iii)$ and $(vi)\Longrightarrow (iii)$ are obvious. The implication $(ii)\Longrightarrow(i)$ is immediate by Proposition \ref{31}, while the implication $(iii)\Longrightarrow(i)$ and the equivalence of $(iii)$ and $(iv)$ follow directly by Proposition \ref{bolg} and Theorem \ref{21}, respectively, since for $(\vY,T)=(\R,\mf{B})$ every measure $P_{W_n}$ on $\mf{B}$ is perfect and $\mf{B}$ is countably generated. The latter together with implication $(vi)\Longrightarrow(iii)$ yields $(vi)\Longrightarrow(iv)$. 

Ad $(ii)\Longrightarrow(iv)$: If $(ii)$ holds true, then there exists a $P_{\vT}$-null set $H_{*}\in{\mf{B}_d}$ such that for any $\theta\notin{{H}}_{*}$ the process $\cip$ is $P_{\theta}$-exchangeable, 
implying together with property (d2) its $P$-exchangeability as well; hence $(iii)$ or equivalently $(iv)$ follows.

Ad $(ii)\Longrightarrow(v)$: Assume that $(ii)$ is true. Putting $S_{\omega}(E):=P_{\theta}(E)$ for any $\omega\in\vO$, $E\in\vS$ and $\theta=\vT(\omega)$, we clearly get that $\{S_{\omega}\}_{\omega\in\vO}$ is a subfield r.c.p. for $P$ over $R:=P\mid\sigma(\vT)$ such that the interarrival times $W_n$, $n\in\N$, are $S_{\omega}$-i.i.d. with a common probability distribution $\mathbf{K}(\omega)$ for any $\omega\notin{H_{**}}:=\vT^{-1}({H_{*}})$, where $\mathbf{K}(\omega):=\mathbf{K}(\theta)$ for each $\omega\in\vO$ and $\vT(\omega)=\theta\notin{H_{*}}\in{\mf{B}_d}$. Since clearly ${H_{**}}$ is an $R$-null set, it follows that $\{S_{\omega}\}_{\omega\in\vO}$, $\mathcal{F}:=\sigma(\vT)$ and $\{\mathbf{K}(\omega)\}_{\omega\in\vO}$ satisfy assertion $(v)$.
\smallskip

Ad $(v)\Longrightarrow(vi)$: 
Assume that $(v)$ holds true and let ${\mathcal F}$, $\{S_{\omega}\}_{\omega\in\vO}$ and $R$ be as in $(v)$. Put ${\wt{\vY}}:=\vO$, $\{S_{\wt{y}}\}_{\wt{y}\in{\wt{\vY}}}:=\{S_{\omega}\}_{\omega\in\vO}$ and $B({\wt{\vY}}):=\sigma(\{S_{\bdot}(E):E\in\vS\})$. Then $B({\wt{\vY}})\subseteq{\mathcal F}$ and so we may define the probability measure $\nu:=R\mid B({\wt{\vY}})$. Since by $(v)$ the process $\{W_n\}_{n\in\N}$ is $S_{\omega}$-i.i.d. for $R$-a.a. $\omega\in\vO$, we get that it is $S_{\wt{y}}$-i.i.d. for $\nu$-a.a. $\wt{y}\in{\wt{\vY}}$. 
The latter together with an application of (sf2) yields that $\cnp$ is a $\nu$-MRP associated with $\{S_{\wt{y}}\}_{\wt{y}\in\vY}$; hence assertion $(vi)$ follows.

Moreover, if $\vS$ is countably generated and $P$ is perfect, the implication 
$(iv)\Longrightarrow(v)$ holds true. In fact, by Theorem \ref{21} we obtain that assertion $(iv)$ is equivalent with the fact that $\{W_n\}_{n\in\N}$ is $P$-conditionally i.i.d. over $\mathcal{F}$. But note that according to Remark \ref{mag}, (a) there exists a subfield r.c.p. $\{S_{\omega}\}_{\omega\in\vO}$ for $P$ over $R:=P\mid\mathcal{F}$. Thus, we may apply Corollary \ref{25p} to get $(v)$.

Assuming now that $(i)$ holds true, it follows by Remark \ref{mag}, (a) that there exists a disintegration
$\{P_{\theta}\}_{\theta\in{\R^{d}}}$ of $P$ over $P_{\vT}$ consistent with $\vT$. So according to Proposition 
\ref{31} 
we get that for $P_{\vT}$-a.a. $\theta\in{\R^d}$ the family $\cnp$ is a $(P_{\theta},\mathbf{K}(\theta))$-RP; hence $(i)$ implies $(ii)$.

Thus, assuming that $\vS$ is countably generated and $P$ is perfect, we obtain that items $(i)$ to $(vi)$ are all equivalent. This completes the whole proof.\hfill$\Box$
\medskip

Note that the most important applications in Probability Theory are still rooted in the case of standard Borel spaces; hence of spaces satisfying always the assumptions of 
the above theorem concerning $P$ and $\vS$. 
\bigskip

\section{Examples}\label{con}

In this section, we provide two groups of examples: The first one (Examples \ref{exa} to \ref{exd}) shows that the perfectness assumption for the probability measures as well as the countability assumption for the $\sigma$-algebras involved in Theorems \ref{21}, \ref{9} and in Corollary \ref{21a}, are essential for the validity of the equivalences obtained therein.

In the second group of examples (Examples \ref{12} to \ref{16}), 
the existence of non-trivial probability spaces admitting $(P,\mathbf{K}(\vT))$-MRPs with prescribed distributions for their interarrival processes as well as for the parameter $\vT$ is 
proven, providing in this way a method of constructing such processes. As a consequence, concrete examples of MRPs are presented. It is also worth noticing that our construction relies on Proposition \ref{31} and allows us to 
obtain probability spaces that satisfy all assumptions of Theorems \ref{21} and \ref{9}. 

It follows an example to show that the perfectness assumption for the probability measures $P_{X_i}$ in Theorem \ref{21} and its Corollary \ref{21a} is essential for the validity of the equivalences of $(i)$ to $(v)$.

\begin{ex}\label{exa}
\normalfont
Let $\vY:=\R_+$ and let $\wt{Q}$ be a probability measure on $\mf{B}(\vY)$. 
Consider a subset $B$ of $\vY$ such that $\wt{Q}^{*}(B)=\wt{Q}^{*}(B^c)=1$, where $\wt{Q}^{*}$ is the outer measure induced by $\wt{Q}$. 
Let $T:=\sigma(\{\mf{B}(\vY),B\})$ and let $Q:T\longrightarrow[0,1]$ be the probability measure defined by
$$
Q(A):=\wt{Q}^{*}(A\cap B)\quad\mbox{for each $A\in T$}.
$$
Then $Q$ is non perfect. Set $\vO:=\vY^{\N}$, $\vS:=T_{\N}$, $P:=Q^{\N}$, $W_n:=\pi_n:\vO\longrightarrow\vY$, where $\pi_n$ is the canonical projection for any $n\in\N$. 
Clearly, assertion $(i)$ of Theorem \ref{21} is satisfied by $\cip$; hence by Theorem \ref{21} we equivalently get that $(ii)$ holds true. 
\smallskip

Assume that assertion $(iii)$ of Theorem \ref{21} is valid, i.e. that there exists a $\sigma$-subalgebra $\mathcal{F}$ of $\vS$ and a family $\{Q_{\omega}\}_{\omega\in\vO}$ of $T$-$\mathcal{F}$-Markov kernels such that
$$
\int_{F}Q_{\omega}^{\N}(H)R(d\omega)=P(F\cap W^{-1}(H))
$$
for every $H\in T_{\N}$ and $F\in\mathcal{F}$, where $R:=P\mid\mathcal{F}$ and $Q_{\omega}^{\N}$, $W$ are as in Theorem \ref{21}. Then $\{Q_{\omega}^{\N}\}_{\omega\in\vO}$ is a subfield r.c.p. for $P$ over $R$; 
hence we may and do assume that $\mathcal{F}$ is countably generated. Applying now a monotone class argument 
we deduce that there exists an $R$-null set $N\in{\mathcal{F}}$ such that for each $A\in\mathcal{F}$ condition $Q_{\omega}^{\N}(A)=1$ holds true for any $\omega\in{N^c}\cap{A}$. 

But by (sf2) we get for every $F\in\mathcal{F}$ that
$$
\int_{F}Q_{\omega}^{\N}(B^{\N})dR=P(F\cap B^{\N})=P(F)=\int_F\chi_FdR,
$$
implying that $P(D)=0$, where $D:=\{\omega\in\vO:Q_{\omega}^{\N}(B^{\N})\neq1\}$.

Put $E:=D\cup{N}$. For any $\omega\in{E^c}$ we get $Q_{\omega}^{\N}(\{\omega\})=1$ and $Q_{\omega}^{\N}(B^{\N})=1$; hence $Q_{\omega}^{\N}(B^{\N}\cap\{\omega\})=1$, implying $B^{\N}\cap\{\omega\}\neq\emptyset$ or $\omega\in{B^{\N}}$. Thus, we get $E^{c}\subseteq B^{\N}$. But then $\wt{Q}^{*}(B^{c})=1$ yields $P(E^c)=0$ or equivalently $P(E)=1$, a contradiction.

Assume now that assertion $(v)$ of Corollary \ref{21a} is valid. It then follows that $P$ must be perfect (see \cite{fa}, Theorem $4^{\prime}$); hence $P_{W_n}$ must do so, a contradiction.
\end{ex}

It follows an example to show that the countability assumption for $T$ or $\vS$ in Corollary \ref{21a} is essential for the validity of the equivalences of $(i)$ to $(v)$ in Theorem \ref{21} and Corollary \ref{21a}.
 
\begin{ex}\label{exb}
\normalfont
Let $\vY:=\R_+$, let $\wt{Q}$ be a probability measure on $\mf{B}(\vY)$, let $T$ be the completion of $\mf{B}(\vY)$ with respect to $\wt{Q}$, and let $Q$ be the completion of $\wt{Q}$. Put $\vO:=\vY^{\N}$, $\wt\vS:=T_{\N}$ and $\wt{P}:=Q^{\N}$. Denote by $\vS$ the completion of $\wt\vS$ with respect to $\wt{P}$ and by $P$ the completion of $\wt{P}$. Then $P$ is a perfect probability measure on $\vS$ (cf. e.g. \cite{fr4}, Proposition 451G and Theorem 451J) but $\vS$ and $T$ are not countably generated. 

Put $W_n:=\pi_n:\vO\longrightarrow\vY$ for any $n\in\N$. Clearly, $\cip$ is $P$-exchangable, that is, assertion $(i)$ of Theorem \ref{21} is valid; hence assertion $(ii)$ of the same theorem is also valid.

Assume now that assertion $(v)$ of Corollary \ref{21a} is valid. It then follows that there exists a $d$-dimensional random vector $\vT$ and a disintegration $\{P_{\theta}\}_{\theta\in\R^d}$ of $P$ over $P_{\vT}$ consistent with $\vT$, implying the existence of a subfiled r.c.p. $\{R_{\omega}\}_{\omega\in\vO}$ of $P$ over $P\mid\sigma(\vT)$, where each function $R_{\bdot}(A):\vO\longrightarrow[0,1]$ is defined by $R_{\omega}(A):=(P_{\bdot}(A)\circ\vT)(\omega)$ for any fixed $A\in\vS$. Since $\sigma(\vT)$ is countably generated, it follows as in Example \ref{exa} that there exists a set $N\in\sigma(\vT)$ such that $P(N)=0$ and for any $A\in\sigma(\vT)$ condition $R_{\omega}(A)=1$ holds true for any $\omega\in{N^c}\cap{A}$.

Choose a set $D\subseteq N^c$ such that $D\notin\sigma(\vT)$ but $D\in\vS$. 
Then for each $\omega\notin{N}$ we obtain
$$
1=R_{\omega}(\{\omega\})\leq R_{\omega}(D)\leq 1\quad\mbox{if}\;\;\omega\in{D}
$$ 
and
$$
1=R_{\omega}(\{\omega\})\leq R_{\omega}(D^c)\leq 1\quad\mbox{if}\;\;\omega\in{D^c}.
$$ 
Thus, $D=N^c\cap\{\omega\in\vO:R_{\omega}(D)=1\}\in\sigma(\vT)$, which is impossible by the choice of $D$; hence assertion $(v)$ of Corollary \ref{21a} is not valid.
\end{ex}

The next example shows that there exists a non perfect probability space $(\vO,\vS,P)$ and a counting process $\cnp$ on it satisfying conditions $(i)$, $(iii)$, $(iv)$ and $(vi)$ of Theorem \ref{9} but not $(ii)$ and $(v)$; hence the perfectness of $P$ is an essential assumption.

\begin{ex}\label{exc}
\normalfont
Let $\vY:=\R_+$, let $Q$ be a probability measure on $T:=\mf{B}(\vY)$, 
and let $(\vO,\wt\vS,\wt{P}):=(\vY^{\N},T_{\N},Q^{\N})$. Consider a subset $B$ of $\vY^{\N}$ such that $\wt{P}^{*}(B)=\wt{P}^{*}(B^c)=1$, the $\sigma$-algebra $\vS:=\sigma(\{\wt\vS\cup B\})$ and the non perfect probability measure $P$ on $\vS$ defined 
by means of 
$$
P(A):=\wt{P}^{*}(A\cap B)\quad\mbox{for each}\quad A\in\vS.
$$
Let $W_n:=\pi_n:\vO\longrightarrow\vY$ be the canonical projection for any $n\in\N$. 
Clearly, assertion $(iii)$ of Theorem \ref{9} is satisfied by $\cip$, and so we equivalently get that assertions $(iv)$ and $(i)$ of the same theorem also hold true. 
Furthermore, it can be easily seen that $(iv)\Longrightarrow (vi)$; hence $(iv)\Longleftrightarrow (vi)$.

Applying similar arguments with those in Example \ref{exa}, we obtain that assertions $(v)$ and $(ii)$ of Theorem \ref{9} are not valid.
\end{ex}
The next example shows that the countability assumption for $\vS$ in Theorem \ref{9} 
is essential for the equivalence of items $(i)$ to $(vi)$ of Theorem \ref{9}.

\begin{ex}\label{exd}
\normalfont
Let $(\vY,T,Q)$ and $(\vO,\wt\vS,\wt{P})$ be as in Example \ref{exc}, let $\vS$ be the completion of $\wt\vS$ with respect to $\wt{P}$ and $P$ the completion of $\wt{P}$, and let $\cip$ be as in Example \ref{exc}. Then $P$ is a perfect probability measure on $\vS$ but $\vS$ is not countably generated. 

It then follows that assertion $(iii)$ of Theorem \ref{9} holds true, while by the same theorem we get that its assertions $(i)$, $(iii)$ and $(iv)$ are all equivalent. Furthermore, it can be easily shown that $(iv)$ implies $(vi)$.
But using similar arguments as in Example \ref{exb} we get that assertion $(ii)$ of Theorem \ref{9} is not valid.

Assume that assertion $(v)$ of Theorem \ref{9} holds true. It then follows that there exist a $\sigma$-subalgebra $\mathcal{F}$ of $\vS$ and a subfield r.c.p. 
$\{S_{\omega}\}_{\omega\in\vO}$ of $P$ over $R:=P\mid\mathcal{F}$.

There exists a countably generated $\sigma$-subalgebra $\mathcal{A}$ of $\mathcal{F}$ such that $\{S_{\omega}\}_{\omega\in\vO}$ satisfies (sf2) with $\mathcal{A}$ in the place of $\mathcal{F}$, and 
\begin{description}
\item[(sf1$^{\prime}$)]
for every $H\in\vS$ there exists a $P$-null set $N_H\in\mathcal{F}$ such that the function 
$S_{\bullet}(H)\mid{(N_H)^c}$ is $\mathcal{A}$-measurable.
\end{description}

In fact, applying similar arguments with those 
in the proof of Proposition \ref{bolg} we get a countably generated $\sigma$-subalgebra $\wt{\mathcal{A}}$ of $\mathcal{F}$ such that for any fixed $\wt{H}\in\wt\vS$ the function {$S_{\bullet}(\wt{H})$} is $\wt{\mathcal{A}}$-measurable. Since for $R$-a.a. $\omega\in\vO$ the measures $P$ and $S_{\omega}$ have the same null
sets, we obtain that for any $H\in\vS$ there exist sets $\wt{H}\in\wt\vS$, $M_H\in\vS$ and $N_H\in\mathcal{F}\cap\vS_0$ such that $H=\wt{H}\cup{M_H}$ and $S_{\omega}(M_H)=P(M_H)=0$ for all $\omega\notin{N_H}$; hence $\{S_{\omega}\}_{\omega\in\vO}$ satisfies (sf1$^{\prime}$) with $\wt{\mathcal{A}}$ in the place of $\mathcal{A}$. 
Applying a monotone class argument we get a $P$-null set $N_1\in\mathcal{F}$ such that for any fixed ${\wt{H}}\in\wt{\mathcal{A}}$ the function ${S_{\bullet}(\wt{H})}\mid{(N_1)^c}$ is $\wt{\mathcal{A}}$-measurable.

Put $\mathcal{A}:=\sigma(\wt{\mathcal{A}}\cup\{N_1\})\subseteq\mathcal{F}$. 
It then follows that $\mathcal{A}$ is
countably generated and $\{S_{\omega}\}_{\omega\in\vO}$ satisfies (sf2) with $\mathcal{A}$ in the place 
of $\mathcal{F}$, as well as (sf1$^{\prime}$). In particular, for any fixed $H\in\mathcal{A}$ the function $S_{\bullet}(H)\mid(N_1)^c$ is $\mathcal{A}$-measurable.  
Applying now a monotone class argument we obtain a $P$-null set $N_2\in\mathcal{A}$ such that for any $H\in\mathcal{A}$ 
and $\omega\in{N^c}\cap{H}$, where $N:=N_1\cup{N}_2$, condition ${S_{\omega}}(H)=1$ holds true.

Choose a set $D\subseteq N^c$ such that $D\notin\mathcal{A}$ but $D\in\vS$.
Then following the same reasoning as in Example \ref{exb}, we deduce that $D\in\mathcal{A}$, which is impossible by the choice of $D$; hence assertion $(v)$ of Theorem \ref{9} is not valid.
\end{ex}
{\em Throughout what follows, we put $\wt\vO:={\R}^{\N}$, $\vO:=\wt\vO\times{\R^{d}}$, 
$\wt\vS:={\mf B}(\wt\vO)$ and $\vS:=\wt\vS\otimes{Z}$ for simplicity}.
The next result extends Theorem 3.1 from \cite{lm5jmaslo}, which provides a construction for MPPs.

\begin{ex}\label{12}
\normalfont
Following the reasoning of Theorem 3.1 from \cite{lm5jmaslo} 
but with $\R^d$, $\mf{B}$ and $Q_n(\theta)=\mathbf{K}(\theta)$ in the place of $\vY:=(0,\infty)$, $\mf{B}(\vY)$ and $Q_n(\theta)=\mathbf{Exp}(\theta)$, respectively, the existence of $(P,\mathbf{K}(\vT))$-MRPs with prescribed distributions for their interarrival processes as well as for the parameter $\vT$ is proven.

In fact, fix on arbitrary $\theta\in{\R^{d}}$. 
If $Q_n(\theta)=\mathbf{K}(\theta)$ for each $n\in\N$, it follows that there exist a unique probability measure $\wt{P}_{\theta}:=\otimes_{n\in\N}Q_n(\theta)$ on $\wt\vS$, and a sequence $\{\wt{W}_n\}_{n\in\N}$ of $\wt{P}_{\theta}$-independent random variables on $(\wt\vO,\wt\vS)$ such that
$$
\wt{W}_n(\omega)=\omega_n=\wt{\pi}_n(\omega)\quad
\mbox{for each}\quad\omega\in\wt\vO\quad\mbox{and}\quad n\in\N,
$$
where $\wt{\pi}_n$ is the canonical projection from $\R^{\N}$ onto $\R$, 
satisfying
$$
(\wt{P}_{\theta})_{\wt{W}_n}=Q_n(\theta)\;\;\mbox{for all}\;\;n\in\N.	
$$
Since by assumption, for any fixed $B\in\mf{B}$ each function $Q_n(\cdot)(B)$ is $Z$-measurable, it follows by a monotone class argument that the same holds true for the function $\wt{P}_{\bdot}(E)$ for fixed $E\in\wt\vS$. 

For each $\theta\in{\R^{d}}$ put $P_{\theta}:=\wt{P}_{\theta}\otimes\delta_{\theta}$,
where $\delta_{\theta}$ is the Dirac probability measure on $Z$, and for each $n\in\N$ set $W_n:=\wt{W}_n\circ{\pi}_{\wt\vO}=\pi_n$, where $\pi_{\wt\vO}$ and $\pi_n$ are the canonical projections from $\vO$ onto $\wt\vO$ and from $\vO$ onto $\R$, respectively. Put now
$$
P(E):=\int\wt{P}_{\theta}(E^{\theta})\mu(d\theta)\quad\mbox{for each}\quad E\in\vS,
$$
where $E^{\theta}:=\{\omega\in\vO:(\omega,\theta)\in E\}$. Then $P$ is a probability measure on $\vS$ such that $\{P_{\theta}\}_{\theta\in{\R^{d}}}$ is a disintegration of $P$ 
over $\mu$ consistent with $\pi_{{\R^{d}}}$, where $\pi_{{\R^{d}}}$ is the canonical projection from $\vO$ onto ${\R^{d}}$ (see \cite{lm5jmaslo}, proof of Theorem 3.1). Furthermore, it can be proven that for all $\theta\in{\R^{d}}$ the process $\{W_n\}_{n\in\N}$ is $P_{\theta}$-independent and $(P_{\theta})_{W_n}=\mathbf{K}(\theta)$ for each $n\in\N$. Clearly, putting $\vT:=\pi_{{\R^{d}}}$ we get $P_{\vT}=\mu$.

It then follows that the \cnpw $\cnp$ induced by $\cip$ is a $(P_{\theta},\mathbf{K}(\theta))$-RP 
for all $\theta\in{\R^{d}}$; hence by Proposition \ref{31} it is a $(P,\mathbf{K}(\vT))$-MRP.
\end{ex}

Applying now Example \ref{12}, we compute the corresponding disintegrating probability measures $P_{\theta}$ ($\theta\in\R^d$) as well as the probability measure $P$ for some MRPs of special interest which are not MPPs. To this aim recall that by $\leb_d$ is denoted the restriction of the Lebesgue measure $\lambda_d$ to $\mf{B}_d$, while any restriction of $\leb_d$ to $\mf{B}(A)$, where $A$ is any Borel subset of $\R^d$, will be denoted again by $\leb_d$. In particular, if $d=1$ then $\leb:=\leb_1=\lambda\mid\mf{B}$, where $\lambda$ is the Lebesgue measure on $\R$.
\smallskip

In the next example, a concrete $(P,\mathbf{K}(\vT))$-MRP is constructed for one of the most common choices that can be made for an interarrival time distribution, i.e. $\mathbf{Ga}(\theta_1,\theta_2)$ with $\theta_1>0$ and $\theta_2=1/2\in(0,1)$, see e.g. \cite{gr}, page 95. In fact, this class of distributions is of special interest, since none of its members satisfy Assumption 5.1 from \cite{gr}, proposed by Huang in \cite{hu}, Theorem 3, which is essential in Grandell's study for MRPs (see \cite{gr}, Section 5.3). Moreover, in the same example it is shown that
there are counting processes $\cnp$ being both $(P,\mathbf{K}(\vT))$-MRPs 
and $P_{\vT}\mid B(\R)$-ones, which are not, though, MRPs according to Grandell \cite{gr}, Definition 5.3.

\begin{ex}\label{12n}
\normalfont
Let $Q_n(\theta)=\mathbf{Ga}(\theta,1/2)$ 
for each $n\in\N$ and for any fixed $\theta>0$, and let $\mu=\mathbf{Ga}(\gamma,\alpha)$. 
Then the conclusions of Example \ref{12} are fulfilled for $d=1$; hence $\wt\vO={\R}^{\N}$, $\vO={\R}^{\N}\times\R$, while $\wt\vS$, $\vS$, $\wt{P}$, $P$, $\{\wt{P}_{\theta}\}_{\theta>0}$, $\{P_{\theta}\}_{\theta>0}$ and $\vT$ are as in Example \ref{12}. 

We first compute the probability measures on measurable cylinders.
Let $\wt{\mathcal{C}}$ denote the family of all measurable cylinders 
$\wt B\in{\mf B}(\wt\vO)$, i.e. of all sets $\wt B\subseteq\wt\vO$ expressible as $\prod_{n\in\N}\wt B_n$, where $\wt B_n\in\mf{B}$ for every $n\in\N$, and 
$\wt{L}:=\{n\in\N:\wt B_n\neq{\R}\}$ is finite. Set $\wt C_n=\wt B_n$ for $n\in\wt{L}$. Then $\wt B=\prod_{k\in\wt{L}}\wt C_k\times{\R}^{\N\setminus\wt{L}}$, so we get
\begin{equation}\label{8a}
\wt{P}_{\theta}(\wt B)=(\otimes_{n\in\N}Q_n(\theta))(\wt B)
=\prod_{k\in\wt L}Q_k(\theta)(\wt C_k)
=\sqrt{\frac{\theta}{\pi}}
\prod_{k\in\wt L}\int_{\wt{C}_k}
\omega_k^{-\frac{1}{2}}e^{-\theta\omega_k}\leb(d\omega_k)
\end{equation}
for each $\theta>0$.
Consider now a measurable cylinder 
$\wt B\times E\in\wt{\mathcal{C}}\times\mf{B}{\bigl((0,\infty)\bigr)}$. 
Applying (\ref{8a}), we get
$$
P_{\theta}(\wt B\times E)
=\wt{P}_{\theta}(\wt B)\delta_{\theta}(E)
=\chi_E(\theta)\sqrt{\frac{\theta}{\pi}}\prod_{k\in\wt L}
\int_{\wt{C}_k}\omega_k^{-\frac{1}{2}}e^{-\theta\omega_k}\leb(d\omega_k);
$$
for each $\theta>0$; hence
$$
P(\wt{B}\times E)
=\frac{\gamma^{\alpha}}{\Gamma(\alpha)\sqrt{\pi}}
\int_{E}
\Bigl[\prod_{k\in\wt L}\int_{\wt{C}_k}\omega_k^{-\frac{1}{2}}e^{-\theta(\gamma+\omega_k)}
\leb(d\omega_k)\Bigr]
\theta^{\alpha-\frac{1}{2}}\leb(d\theta).
$$
As a consequence, by applying standard methods of Topological Measure Theory, the probability measures $P(E)$ and $P_{\theta}(E)$ can be computed for any $E\in\vS$. For details see \cite{lm5jmaslo}, Example 3.3, (b).
\end{ex}

Finally, it follows an example to show that we cannot avoid including in Huang's definition of an MRP the assumption that $\cip$ is $P_{\wt{y}}$-identically distributed for $\nu$-a.a. $\wt{y}\in{\wt{\vY}}$ (see also Remark \ref{27}).

\begin{ex}\label{16}
\normalfont
Let $d=1$. If $Q_n(\theta)=\mathbf{Exp}(n\theta)$ for each $n\in\N$ and for any fixed $\theta>0$, 
and if $\mu=\mathbf{Ga}(2,1)$ then all requirements of Example \ref{12} except for 
$$
Q_n(\theta)=\mathbf{K}(\theta)\quad\mbox{for all $n\in\N$ and for any fixed $\theta>0$}
$$ 
are satisfied. In fact, in this case $\mathbf{K}(\theta)$ is substituted by $\mathbf{K}(n\theta):=\mathbf{Exp}(n\theta)$. 

So the probability measures $\wt{P}$ and $P$ on $\wt\vS=\mf{B}(\wt\vO)$ and $\vS=\mf{B}(\vO)$, where $\wt\vO={\R}^{\N}$ and $\vO=\R^{\N}\times{\R}$, respectively, as well as the disintegrations
$\{\wt{P}_{\theta}\}_{\theta>0}$ and $\{P_{\theta}\}_{\theta>0}$ can be computed. Moreover, there exists a random variable $\vT$ on $\vO$ such that $P_{\vT}=\mathbf{Ga}(2,1)$. Following the same reasoning as in Example \ref{12}, we also obtain an interarrival process $\cip$ which is $P_{\theta}$-independent for all $\theta>0$ and satisfies $W_n=\pi_n$ as well as $(P_{\theta})_{W_n}=\mathbf{K}(n\theta)$ for all $n\in\N$ and for any fixed $\theta>0$.

But the $P_{\theta}$-independence of $\cip$, for all $\theta>0$, implies
for every $r\in\N$ and for all $w_1,\ldots,w_r\in\R_+$ that
\begin{equation}\label{15a} 
P\Bigl(\bigcap_{k=1}^{r}\{W_k\leq w_k\}\Bigr)
=\int\prod_{k=1}^{r}P_{\theta}(W_k\leq w_k)\nu(d\theta), 
\end{equation}
where $\nu=P_{\vT}\mid B((0,\infty))=\mu\mid B((0,\infty))$ and $B((0,\infty))=\sigma(\{P_{\bdot}(E):E\in\vS\})$.
So, $\cip$ is an interarrival process which is not $P_{\theta}$-identically distributed for any fixed $\theta>0$ but which satisfies (\ref{15a}). As a consequence, the \cnpw $\cnp$ induced by the sequence of canonical projections $\{\pi_n\}_{n\in\N}=\cip$ is not a $\nu$-MRP associated with $\{P_{\theta}\}_{\theta\in\vX}$. Furthermore, for every $w_1,w_2\in\R_+$ we have 
\begin{eqnarray*}
P(W_1\leq w_1,W_2\leq w_2)&=&2\int_{0}^{\infty}(1-e^{-\theta w_1})(1-e^{-2\theta w_2})
e^{-2\theta}d\theta\\
&=&w_2(w_2+1)^{-1}-2[(w_1+2)^{-1}-(w_1+2w_2+2)^{-1}],\label{15c}
\end{eqnarray*}
implying that $P(W_1\leq2,W_2\leq1)=\frac{1}{3}\neq\frac{2}{7}=P(W_1\leq1,W_2\leq2)$; hence $\cip$ is not $P$-exchangeable. 
\end{ex}


{\small
\noindent{\sc D.P. Lyberopoulos and N.D. Macheras}\\
{\sc Department of Statistics and Insurance Science}\\
{\sc University of Piraeus, 80 Karaoli and Dimitriou street}\\
{\sc 185 34 Piraeus, Greece}\\
E-mail: {\tt dilyber@webmail.unipi.gr}$\;$ {\sc and}$\;$ {\tt macheras@unipi.gr}
}
\end{document}